\documentclass[a4paper,10pt]{article}
\usepackage[utf8]{inputenc}
\usepackage[T1]{fontenc}
\usepackage{amsmath,amsthm,amsfonts,amssymb,euscript,stmaryrd,graphicx,enumitem,yhmath,mathrsfs}
\usepackage[toc,page]{appendix}
\usepackage{algorithm,algorithmic}
\usepackage{fullpage}
\usepackage[usenames,dvipsnames]{color}
\usepackage{color,psfrag}
\usepackage{dsfont}
\usepackage[babel=true]{csquotes} 
\usepackage{lmodern}
\usepackage{bibentry}

\usepackage[french,english]{babel}
\usepackage{ifthen}
\usepackage{array}
\usepackage{hyperref} 
\usepackage{verbatim}
\usepackage[normalem]{ulem}

\usepackage{xcolor}
\definecolor{thmcolor}{rgb}{0.71,0.14,0.07}
\usepackage[colored]{shadethm}

\usepackage{empheq}

\newcommand{\footnoteremember}[2]{
 \footnote{#2}
 \newcounter{#1}
 \setcounter{#1}{\value{footnote}}
}

\author
{        {Joachim \textsc{Lebovits}\footnote{MAThematic Center of Heidelberg, University of Heidelberg, INF 294,
69120 Heidelberg, Germany.}\footnoteremember{footnoteP13}{Laboratoire Analyse, Géométrie et Applications, C.N.R.S. (UMR 7539), Université Paris 13, Sorbonne Paris Cité, 99 avenue Jean-Baptiste Clément
93430, Villetaneuse, France. Email address: \url{jolebovits@gmail.com}}}
  }

\newtheorem{theo}{Theorem}[section]
\newtheorem{theodef}{Theorem-Definition}[section]
\newtheorem{prop}[theo]{Proposition}

\newtheorem{rem}{Remark}
\newtheorem{lem}[theo]{Lemma}
\newtheorem{defi}{Definition}

\newtheorem{ex}[theo]{Example}

\newenvironment{pr}{\begin{proof}[\bfseries \textup{Proof.}]}{\end{proof}}

\newenvironment{g}[2]{\begin{array}{cl} %
<\hspace{-0.2cm}<\hspace{-0.1cm} #1, \hspace{-0.3cm} &#2 \hspace{-0.1cm} >\hspace{-0.2cm}>}%
{\end{array}}

\newcommand{\bD}{\mathbb{D}}

\newcommand{\bN}{\mathbb{N}}

\newcommand{\sifbm}{\mathbf{B}}

\newcommand{\E}{\mathbf{E}}

\newcommand{\N}{\mathbf{N}}

\newcommand{\R}{\mathbf{R}}

\newcommand{\zZ}{\mathbf{Z}}

\newcommand{\oS}{${(\cS)}^*$}
\newcommand{\ooS}{{(\cS)}^*}

\newcommand{\Z}{\mathcal{Z}}

\newcommand{\w}{\Gamma_{\hspace{-0.1cm}\scriptscriptstyle R}}

\newcommand{\zr}{\mathcal{Z}_{\scriptscriptstyle R}}

\newcommand{\wa}{\Gamma^{\scriptscriptstyle a,b}_{\hspace{-0.1cm}\scriptscriptstyle R}}

\newcommand{\z}{\mathcal{Z}^{\scriptscriptstyle \hspace{-0.01cm} T}_{\hspace{-0.05cm}\scriptscriptstyle R}}

\newcommand{\za}{\mathcal{Z}^{\scriptscriptstyle \hspace{-0.01cm} a,b}_{\hspace{-0.05cm}\scriptscriptstyle R}}

\newcommand{\di}{d^{\diamond} \hspace{-0.05cm}}

\newcommand{\ds}{\displaystyle}

\def\i1{\mathbf{1}}  

\newcommand{\cS}{\mathcal{S}}
\newcommand{\sS}{\mathscr{S}}

\def\Aa{\ref{Aa} }
\def\Aap{\ref{Aa}}

\def\Aaiii{\ref{Aiii} }

\def\ie{{\textit{i.e. }}}

\def\rR{\mathscr {R}}

\def\cB{{\cal B}}
\def\cC{C}
\def\cD{{\cal D}}

\def\cG{{\mathcal G}}

\def\cS{{\cal S}}

\def\si{\text{sign}}  

\newcommand{\cF}{\mathcal{F}}

\newcommand{\bit} {\begin{itemize} }
\newcommand{\eit} {\end{itemize} }
\def\ben{\begin{enumerate}}
\def\een{\end{enumerate}}

\def\iti{\item[(i)]}

\def\ita{\item[a)]}

\def\itii{\item[(ii)]}

\def\itb{\item[b)]}

\def\itiii{\item[(iii)]}

\def\ {\hspace{0.1cm}}

\def\keywordname{{\bf Keywords:}}

\newcommand{\keywords}[1]{\par\addvspace\baselineskip\noindent\keywordname\enspace\ignorespaces#1}

{\end{array}}

\makeatletter
\renewcommand\theequation{\thesection.\arabic{equation}}
\@addtoreset{equation}{section}
\makeatother

\graphicspath{{./}{figures/}}

\DeclareFontFamily{U}{mathx}{\hyphenchar\font45}
\DeclareFontShape{U}{mathx}{m}{n}{
      <5> <6> <7> <8> <9> <10>
      <10.95> <12> <14.4> <17.28> <20.74> <24.88>
      mathx10
      }{}
\DeclareSymbolFont{mathx}{U}{mathx}{m}{n}
\DeclareFontSubstitution{U}{mathx}{m}{n}
\DeclareMathAccent{\widecheck}{0}{mathx}{"71}
\DeclareMathAccent{\wideparen}{0}{mathx}{"75}

\title{%
    \begin{minipage}\linewidth
        \centering\bfseries\sffamily
        Local Times of Gaussian Processes
        \vskip3pt
        \large Stochastic Calculus with respect to Gaussian Processes: Part II
            \end{minipage}
}

\begin{document}

\maketitle

\frenchspacing

\vspace{-5ex}

\begin{abstract}
\begin{footnotesize}
\noindent The aim of this work is to define and perform a study of local times of all Gaussian processes that have an integral representation over a real 
interval (that maybe infinite). Very rich, this class of Gaussian processes, 
contains Volterra processes (and thus fractional Brownian motion),  
multifractional Brownian motions as well as processes, the regularity of which varies along the time. 

\noindent Using the White Noise-based anticipative stochastic calculus with respect to Gaussian processes developed in \cite{JL17a}, we first establish a Tanaka formula. This allows us to define both weighted and non-weighted local times and finally to provide occupation time formulas for both these local times. A complete comparison of the Tanaka formula as well as the results on Gaussian local times we present here, is made with the ones proposed in 
 \cite{MV05,LN12,SoVi14}. 

\end{footnotesize}
\end{abstract}

\vspace{-3ex}


\keywords{
\begin{footnotesize}
{\small Gaussian processes and local times, stochastic analysis \& White noise theory, Wick-Itô integrals,  Tanaka formulas, fractional and multifractional Brownian motions.}
\end{footnotesize}
}
\smallskip


{\small \bfseries AMS Subject Classification:} 60G15; 60H40; 60H05; 60G22

\vspace{-2ex}

\section{Introduction}
\vspace{-1ex}
The purpose of this paper is to define and perform a study of local times of all Gaussian processes that have an integral representation over a real 
interval (that maybe infinite). Theis class of Gaussian processes, denoted $\mathscr{G}$, is defined as being the set of Gaussian processes $G:={(G_{t})}_{t\in\rR}$ which can be written under the form:
\vspace{-1ex}
\begin{equation}
\label{erofkorekp}
G_{t}=\int_{\R} g_{t}(u) \ dB_{u} ,
\end{equation}
where $\R$ denotes the set of real numbers, $\rR$ denotes a closed 
interval of $\R$ (that may be equal to $\R$), $B:={(B_{u})}_{u
\hspace{-0.01cm}\in\R}$ is Brownian motion on $\R$ and  ${(g_{t})}_{t\in
\rR}$  is a family of a measurable square integrable functions
on $\R$. 
The set $\mathscr{G}$ obviously contains Volterra processes, as well as Gaussian Fredholm processes\footnote{See \cite[Introduction]{JL17a} for a precise definition of these processes.}. 

In order to perform this study of Gaussian local times we will take advantage of the White Noise-based anticipative stochastic calculus with respect to Gaussian processes in $\mathscr{G}$, developed in \cite{JL17a}.
Our main results are:
\vspace{-1ex}
\bit
\item A Tanaka formula, that reads, for every $T>0$: 
\eit
\vspace{-1ex}
\begin{equation*}
|G_{t} - c| = |c| + \int^{T}_{0} \   \text{sign}(G_{t} - c) \  \di G_{t} +  \int^{T}_{0}  \  \delta_{\{c\}}(G_{t}) \   dR_{t},
\end{equation*}
where the equality holds in $L^2(\Omega)$, where $t\mapsto R_{t}$ denotes the variance function of $G$, which will be supposed to be a continuous function, of bounded variations; the meaning of the different terms will be explained below. 
\vspace{-1ex}
\bit
\item Two occupation time formulas that read, for every positive real-valued Borel function $\Phi$:
\eit
\vspace{-3ex}
\begin{align*}
&  \int^{T}_{0}\  \Phi(G_{s}(\omega)) \ ds =  \int_{\R}\  \ell^{(G)}_{T}(y)(\omega) \cdot \Phi(y) \ dy;\\
\vspace{-0.5cm}  &\int^{T}_{0}\  \Phi(G_{s}(\omega)) \ dR_{s} =  \int_{\R}\  \mathscr{L}^{(G)}_{T}(y)(\omega) \cdot \Phi(y) \ dy,
\end{align*}
where $\ell^{(G)}_{T}(y)$ (\textit{resp.} $\mathscr{L}^{(G)}_{T}(y)$)  denotes the non weighted (\textit{resp.} the weighted) local time of $G$ (which will be both defined precisely), at point $y$,  up to time $T$. We will moreover establish that these two latter random variables are $L^{2}$ random variables.

If this work is self-contained, it would be helpful for many readers to refer to \cite{JL17a} to get a precise idea of what White Noise-based stochastic calculus \textit{wrt} Gaussian processes is.
\smallskip

\textit{Outline of the paper }
The remaining of this paper is organized 
as follows.  In Section 
\ref{qdokpsdksqopdqskopqsdkopdqs}, we recall basic
facts about white noise theory and about the family of 
operators ${(M_{H})}_{H\in(0,1)}$, which is instrumental for 
our running example. Section \ref{Itô} is devoted to the obtention of an Itô formula for generalized functionals of $G$.
This result is a preparatory work for Section \ref{Tanaka}. More precisely, 
a Tanaka formula is given in  Subsection \ref{tanakaformula}, while the definition of weighted and non weighted local 
times of elements of $\mathscr{G}$, as well as their integral representation is given in Subsection \ref{ozdeiojezoijziojéazpokdzlkn}.  An 
occupation time 
formula for both, weighted and non weighted, local 
times of $G$ is established in Subsection \ref{occupation} and constitutes the second main result of this paper. We also prove in Section \ref{occupation} that both, weighted and non weighted local times of $G$, belong, as 
two parameters processes, to $L^{2}(\R\times\Omega, 
\lambda \otimes \mu)$. Finally, a complete comparison of the Tanaka formula as well as the results on Gaussian local times we present here, with the ones proposed in 
 \cite{MV05,LN12,SoVi14} is made in Subesction \ref{tanakaformula} for Tanaka formula and at the end of Subsection \ref{comparaison} for the other results.

%
%

\section{Background on White Noise Theory and White Noise-based stochastic calculus
}
\label{qdokpsdksqopdqskopqsdkopdqs}

%
%
%
 Introduced by T. \hspace{-0.1cm}Hida in \cite{Hi75}, White Noise Theory is, roughly speaking, the stochastic analogous of deterministic generalized functions (also known as tempered distributions). The idea is to realize nonlinear functional on a Hilbert space as functions of white noise (which is defined as being the time derivative of Brownian motion). 
 We recall in this section the minimum standard set-up for classical white-noise theory. 
  Readers interested in more details 
  may refer  to \cite{HKPS,Kuo2} and \cite{Si12}.
%
%
\subsection{The spaces of stochastic test functions and stochastic distributions}
\label{cddcwcdsd}

Define $\N$ (resp. $\N^*$) the set of non negative integers (resp. positive integers). Let  $\sS(\R)$ be the Schwartz space endowed with its usual topology.
%
%
%
 Denote
$\sS'(\R)$ the space of tempered distributions, which is 
the dual space of $\sS(\R)$, and $\widehat{F}$ 
or $\cF(F)$ 
the 
Fourier transform of any element $F$ of $\sS'(\R)$. For every positive 
real $p$, denote $L^{p}(\R)$ the set of measurable functions $f$ such 
that $\int_{\R} {|f(u)|}^{p} \ du<+\infty$. When $f$ belongs to $L^1(\R)$, $
\widehat{f}$ is defined on $\R$ by setting $
\widehat{f}(\xi) := \int_{\R} e^{-ix\xi}  f(x) \ dx$.
Define the measurable space $(\Omega,\cF)$ by setting $\Omega:= {\sS'}(\R)$ and 
$\cF := \cB({\sS'}(\R))$, where $\cB$ denotes the $\sigma$-algebra of Borel 
sets.
The  Bochner-Minlos theorem ensures that there exists a unique probability measure $\mu$ on $(\Omega,\cF)$ such that, 
for every $f$ in $\sS(\R)$, the map $<.,f>:(\Omega,\cF) \rightarrow \R$ defined by $<.,f>(\omega) = <\omega,f>$ (where 
$<\omega,f>$ is by definition   $\omega(f)$, \ie the action of $\omega$ on $f$) is a centred 
Gaussian random variable with variance equal to ${\|f\|}^2_{L^2(\R)}$ under $\mu$. The map $f \mapsto <.,f>$ being an isometry from 
$(\sS(\R),{< , >}_{L^2(\R)})$ to $(L^2(\Omega,\cF,\mu), {< , >}
_{L^2(\Omega,\cF,\mu)})$, it may be extended to $L^2(\R)$. One may thus consider the centred Gaussian random variable $<.,f>$, for any $f$ in $L^2(\R)$.
In particular, let $t$ be in $\R$, the indicator function ${\i1}_{[0,t]}$ is defined by setting: $\i1_{[0,t]}(s):=1$ if $0 \leq s \leq t$, $\i1_{[0,t]}(s):=-1$ if  $t \leq s \leq 0$ and  $ \i1_{[0,t]}(s) :=0$ otherwise.
%
%
%
%
Then the process ${({\widetilde{B}}_t)}_{t \in \R}$, where ${\widetilde{B}}_t (\omega):= {\widetilde{B}}(t,\omega):= {<\omega, \i1_{[0,t]}>}$ is a standard Brownian motion with respect to $\mu$. It then admits a continuous version which will be denoted $B$.
Define, for $f$ in $L^2(\R)$, $I_1(f)(\omega) := {<\omega, f>} $. Then $I_1(f)(\omega)  =\int_{\R} f(s) \ dB_s (\omega) \hspace{0.2cm}  \mu-{\text {a.s.}}$, where  $\int_{\R} f(s)\ dB_s$ denotes the Wiener integral of $f$.
For every $n$ in $\N$, let $e_n(x):=   {(-1)}^n \  {\pi}^{-1/4} {(2^n n!)}^{-1/2} e^{x^2/2} \frac{d^n}{dx^n}( e^{-x^2})$ be the $n$-th 
Hermite function. It is well known (see $\cite{Tha}$) that ${(e_k)}_{k\in \N}$ is a family of functions of $\sS(\R)$ that forms an 
orthonormal basis of $L^2(\R,dt)$.
The following properties about the Hermite functions (the proof of which can be found in \cite{Tha}) 
will be useful.
\begin{theo}
\label{ozdicjdoisoijosidjcosjcsodijqpzeoejcvenvdsdlsiocfuvfsosfd}
There exist positive constants $C$ and $\gamma$ such that,  for every $k$ in $\N$,
\begin{equation*}
 {|e_k(x)|  } \leq C \ \big( {(k+1)}^{-1/12}  \cdot {\i1}_{\{|x| \leq 2\sqrt{k+1}\}} + e^{-\gamma x^2}\cdot {\i1}_{\{|x| > 2\sqrt{k+1}\}}\big).
\end{equation*}
%
%
\end{theo}

 Let  ${({| \ |}_p)}_{p\in\zZ}$ be the family  norms defined by ${|f|}
 ^2_p:=   \sum^{+\infty}_{k=0} {(2k + 2)}^{2p} \  {<f,e_k>}^2_{L^2(\R)}
 $, for all $(p,f)$ in $\zZ \times L^2(\R)$. The operator $A$, defined on $\sS(\R)$, by setting $A:= -\frac{d^2}{dx^2} + 
x^2 +1$, admits the sequence ${(e_n)}_{n \in \N}$ as eigenfunctions 
and the sequence $({{2n+2}})_{n \in \N}$ as  eigenvalues. 
Define, for $p$ in $\N$, the spaces $\sS_p(\R):=\{f \in L^2(\R), \  {|f|}_{p} <+\infty \}$ and  $
\sS_{\hspace{-0.15cm}-p}(\R)$ as being the completion of $L^2(\R)$ with 
respect to the norm ${{|\  \  |}_{-p}}$.
%
%
We summarize here the minimum background on White Noise Theory, written \textit{e.g.} in \cite[p. 692-693]{LLVH}. More precisely, let $(L^2)$ 
denote the space $L^2(\Omega,\cG,\mu)$, where $\cG$ is the $\sigma$-
field generated by ${(<.,f>)}_{f \in L^2(\R)}$. According to  Wiener-Itô's theorem, for every random variable $
\Phi$ in $(L^2)$ there exists a 
unique sequence ${(f_n)}_{n \in \N}$ of functions in ${\widehat{L}}
^2(\R^n)$ such that $\Phi$ can be decomposed as $\Phi =  {\sum^{+ \infty}
_{n = 0} I_n(f_n)}$, where ${\widehat{L}}^2(\R^n)$ denotes the set of all 
symmetric functions $f$ in $L^2(\R^n)$ and $I_n(f)$ denotes the $n-$th 
multiple Wiener-Itô integral of $f$ with the convention that $I_0(f_0) = f_0$ 
for constants $f_0$.  
For any $\Phi:=  {\sum^{+ 
\infty}_{n = 0}  \ I_n(f_n})$ satisfying the condition ${\sum^{+ \infty}_{n = 0}  
n!\ {|A^{\otimes n}f_n|}^2_{0} }< +\infty$, define the element $\Gamma(A)
(\Phi)$ of $(L^2)$ by $\Gamma(A)(\Phi):= {\sum^{+ \infty}_{n = 0} \ 
I_n(A^{\otimes n} f_n)}$, where $A^{\otimes n}$ denotes the $n-$th tensor 
power of the operator $A$ (see \cite[Appendix E]{Jan97} for more details 
about tensor products of operators).
The operator $\Gamma(A)$ is densely defined on $(L^2)$. It is invertible 
and its inverse  ${\Gamma(A)}^{-1}$ is bounded.
We note, for $\varphi$ in $(L^2)$, ${\|\varphi\|}^2_0:={\|\varphi\|}^2_{(L^2)}$. 
For $n$ in $\N$, let $\bD\text{om}({\Gamma(A)}^n)$ be the domain of the  
$n$-th iteration of $\Gamma(A)$. Define the  family  of norms  ${({\|\ \|}_p)}
_{p \in \zZ}$ by:
\begin{equation*}
 {\|\Phi\|}_p :=  {\|\Gamma(A)^p\Phi\|}_{0} =  {\| \Gamma(A)^p(\Phi)\|}_{(L^2)},  \hspace{1cm} \forall p \in \zZ,\hspace{0.5cm} \forall\Phi \in (L^2)\cap \bD{\text{om}}({\Gamma(A)}^p).
\end{equation*}
For $p$ in $\N$, define $({\cS}_{p}):=\{\Phi \in (L^2): \   \Gamma(A)^p(\Phi) \  \text{exists and belongs to}  \  (L^2) \}$ and define $({\cS}_{-p})$ as the completion of the space  $(L^2)$ with respect to the norm ${{\|\ \|}_{-p}}$.
As in \cite{Kuo2}, we let $(\cS)$ denote the projective limit of the sequence ${( (\cS_{p}))}_{p \in \N}$ and ${(\cS)}^*$ the inductive limit of the sequence  ${(({\cS_{-p}}))}_{p \in \N}$. This means in particular that ${(\cS)} \subset (L)^{2} \subset {(\cS)}^*$ and that ${(\cS)}^*$ is the dual space of $(\cS)$. Moreover, 
$(\cS)$ is called the space of stochastic test functions while ${(\cS)}^*$ the Hida distribution space. We will note  $<\hspace{-0.2cm}<\hspace{-0.1cm} \ ,\  \hspace{-0.1cm}>\hspace{-0.2cm}>$  the 
duality bracket between ${(\cS)}^*$ and $(\cS)$. If $\phi, \Phi$ belong to $(L^2)$, then we have the equality $<\hspace{-0.2cm}<\hspace{-0.1cm} \ \Phi , \varphi\  \hspace{-0.1cm}>\hspace{-0.2cm}> = {<\Phi,\varphi>}_{(L^2)} = \E[\Phi \  \varphi]$. Besides, denote $< , >$  the duality bracket between $\sS'(\R)$ and $\sS(\R)$ and recall that every tempered distribution $F$ can be written as $F = {\sum^{+ \infty}_{n = 0} \ <F,e_n>}\ e_n$, where the convergence holds in $\sS'(\R)$. The next proposition, that will be used extensively in the sequel, is a consequence of the definition of $(\cS)$ and  ${(\cS)}^*$.
\begin{prop}
\label{laiusdhv}
Let $F$ be in in $\sS'(\R)$. Define  $<.,F>:={\sum^{+ \infty}_{n = 0} \  
<F,e_n>\  <.,e_n>}$. Then there exists $p_0$ in $\N$ such that  that $<.,F>
$ belongs to  $({\cS_{-p_0}})$,  and hence to ${(\cS)}^*$. Moreover we have 
${\|<.,F>\|}^2_{-p_0}  = {|F|}^2_{-p_0}$.
    Conversely, define $\Phi:={\sum^{+ \infty}_{n = 0} \ b_n <.,e_n>}$, where $
{(b_n)}_{n \in \N}$ belongs to ${\R}^{\N}$. Then $\Phi$ belongs to ${(\cS)}^*
$ if and only if there exists an integer $p_0$ in $\N$ such that ${\sum^{+ 
\infty}_{n = 0}  \ b^2_n \ {(2n+2)}^{-2p_0}}  < +\infty$. In this latter case $F:= 
{\sum^{+ \infty}_{n = 0} \ b_n e_n}$   belongs to  $\sS_{\hspace{-0.15cm}-
p_0}(\R)$ and then to $\sS'(\R)$. It moreover verifies the equality ${|F|}^2_{-p_0} = 
{\sum^{+ \infty}_{n = 0}  \ {b^2_n} {(2n+2)}^{-2p_0} } = {\|\Phi\|}^2_{-p_0} $.

\end{prop}

\subsection{\oS-process, \oS-derivative and \oS-integral}

%

Let  $(\R,\cB(\R),m)$ be a sigma-finite measure space. Through this section, $I$ denotes an element of $\cB(\R)$.
 A measurable function $\Phi:I \rightarrow $\oS$ $ is called a stochastic distribution process, or an \oS-process. An  \oS-process $\Phi$ is 
said to be differentiable at $t_0 \in I$ if {$\lim\limits_{r 
\to 0} \ r^{-1}\ (\Phi_{t_0+r}-\Phi_{t_0})$} exists in \oS. We 
note $\frac{d\Phi_{t_{0}}}{dt}$ the \oS\hspace{-0.1cm}-
derivative at $t_0$ of the stochastic distribution process $
\Phi$. 
%
It is also possible to define an \oS-valued integral in the following way (one may refer to \cite[p.$245$-$246$]{Kuo2}  or \cite[Def. $3.7.1$ p.$77$]{HP} for more details). 
\begin{theodef}[{\bfseries integral in  \oS}]
\label{izeufheriduscheirucheridu}
Assume that $\Phi:I \rightarrow \ooS$ is weakly in $L^{1}(I,m)$, \ie assume that for all $\varphi$ in  ${(\cS)}$, the 
mapping $u \mapsto\ <\hspace{-0.2cm}<\hspace{-0.1cm} \  \Phi_{u},\ 
\varphi \hspace{-0.1cm}>\hspace{-0.2cm}>$,  from $I$ to $\R$, belongs to  
$L^{1}(I,m)$. Then there exists an unique element in \oS, noted $\int_{I} 
\Phi_{u} \ m(du)$, 
such that,  for all  $\varphi$ in $(\cS)$,
\begin{equation*}
{<\hspace{-0.2cm}<\hspace{-0.1cm} \ \int_{I} \Phi(u) \ m(du), \varphi \ \hspace{-0.1cm}>\hspace{-0.2cm}>} = 
\int_{I} <\hspace{-0.2cm}<\hspace{-0.1cm} \ \Phi_{u},
\varphi\  \hspace{-0.1cm}>\hspace{-0.2cm}> \ m(du). 
\end{equation*}
 \end{theodef}

%

We say in this case that $\Phi$ is \oS-integrable on $I$ (with respect to the measure $m$), in the {\it Pettis sense}. In the sequel, unless otherwise specified, 
we will always refer to the integral in Pettis' sense

\subsection{S-transform and Wick product}
\label{ksksksks}

For $f$ in $L^2(\R)$, define the \textit{Wick exponential} of $<.,f>$, noted $:e^{<.,f>}:$, as the $(L^2)$ random variable equal to
$e^{<.,f> - \frac{1}{2}{|f|}^2_0}$. 
The $S$-transform of an element $\Phi$ of $(\cS^*)$, noted $S(\Phi)$, is defined as the function from  $\sS(\R)$ to $\R$ given by $S(\Phi)(\eta):= {\begin{g}{\Phi}{\ \hspace{-0.35cm}:e^{<.,\eta>}:\hspace{0.075cm}}  \end{g} }$ for any $\eta$ in $\sS(\R)$.
For any $(\Phi,\Psi)\in \ooS\times \ooS$, there exists a unique element of \oS, called the Wick product of  $\Phi$ and $\Psi$,  and  noted $\Phi \diamond \Psi$,  such that $S(\Phi\diamond \Psi)(\eta) = S(\Phi)(\eta) \  S(\Psi)(\eta)$ for every $\eta$ in $\sS(\R)$.  
Note that, when $\Phi $ belongs to $(L^2)$, $S\Phi(\eta)$ is nothing but $\E[\Phi :e^{<.,\eta>}: ] =  e^{ -\frac{1}{2}{|\eta|}^2_0} \  \E[\Phi \  e^{<.,\eta>} ]$. The following result
 will be intensively used in the sequel.
\vspace{-0.1cm}

\begin{lem}{\cite[Lemma 2.3.]{JLJLV1}}
\label{dede}
The inequality $|S(X \diamond Y)(\eta)| \leq {\|X\|}_{-p} \  {\|Y\|}_{-q} \  e^{{|\eta|}^2_{\max \{p;q\}}}$ is valid for any $(p,q)$ in $\N^2$ and $(X,Y)$ in $({\cS}_{-p}) \times ({\cS}_{-q})$.

\end{lem}
%
%




Some useful properties of S transforms are listed in the proposition below. The proof of the results stated in this proposition can be found in \cite[Chap $5$]{Kuo2}.

\begin{prop}[Some properties of S transforms]
\label{qmqmqmmmmm}
When $\Phi$ is deterministic then $\Phi \diamond 
\Psi = \Phi \ \Psi$, for all $\Psi$ in \oS.  Moreover, 
 let $\Phi=\sum^{+\infty}_{k = 0}  a_k\hspace{-0.1cm} <\hspace{-0.1cm}.,\ \hspace{-0.1cm}e_k>$ and $\Psi=\sum^{+\infty}_{n = 0}  I_{n}(f_{n})$ be in  \hspace{-0.1cm} \oS. \hspace{-0.2cm} Then their S-transform is given, for every $\eta$ in $\sS(\R)$, by $S(\Phi)(\eta) = \sum^{+\infty}_{k = 0} a_k\ {<\eta,e_k>}_{L^2(\R)}$\  and $S(\Psi)(\eta) = \sum^{+\infty}_{k = 0} \ <f_{n},\eta^{\otimes n}>$.
Finally, for every $(f,\eta,\xi)$ in $L^{2}(\R) \times \sS(\R)\times \R$, we have the equality:
\begin{equation}
\label{oeij12}
S(e^{i\xi <.,f>})(\eta) = e^{\frac{1}{2}({|\eta|}^{2}_{0} + 2 i \xi <f,\eta> -{\xi}^{2} {|f|}^{2}_{0})}.
\end{equation}
\end{prop}

One may refer to \cite[Chap.$3$ and $16$]{Jan97}  for more details about Wick product. The following results on the S-transform will be used extensively in the sequel. See \cite[p.$39$]{Kuo2}  and \cite[p.280-281]{HKPS}  for proofs. Denote $\cF(A;B)$ 
the set of $B$-valued functions defined on $A$.

\begin{lem}
\label{dkdskcsdckksdksdmksdmlkskdm} 
The $S$-transform verifies the following properties:
\bit
\iti The map $S:\Phi\mapsto S(\Phi)$, from \oS into $\cF(\sS(\R);\R)$, is injective.
\itii Let $\Phi:I\ \rightarrow \ooS$ be an \oS process. If  $\Phi$ is \oS-integrable over $I$ \textit{wrt} $m$, then one has,  for all $\eta$ in $\sS(\R)$,
 $S(\int_{I} \Phi(u) \  m(du) )(\eta) = \int_{I} S(\Phi(u)) (\eta) \  m(du)$.
\itiii  Let $\Phi:I \rightarrow \ooS$  be an \oS-process differentiable at $t\in I$. Then, for every $\eta$ in $\sS(\R)$ the map $u\mapsto [S \Phi(u)](\eta)$ is differentiable at $t$ and verifies $\displaystyle{S[\tfrac{d\Phi}{dt}(t)](\eta) = \tfrac{d}{dt}\big[ S[ \Phi(t)](\eta) \big]}$.
\eit
\end{lem}
The next theorems provide a criterion for integrability in \oS,  in term of $S$-transform. 

\begin{theo}{\cite[Theorem 13.5]{Kuo2}}
\label{peodcpdsokcpodfckposkcdpqkoq}
Let $\Phi:I \rightarrow {(\cS)}^*$ be a stochastic 
distribution such that, for all $\eta$ in $\sS(\R)$, the real-
valued map $t\mapsto S[\Phi(t)](\eta)$ is measurable and such
that there exist a natural integer $p$, a real $a$ and a  
function $L$ in $L^1(I,m)$ such that $ |S(\Phi(t))(\eta)| 
\leq L(t) \  e^{a {|\eta |}^2_{p} }$, for all $\eta$ of $\sS(\R)
$ and for almost every $t$ of $I$. Then $\Phi$ is \oS-
integrable over $I$, \textit{wrt} to $m$.
\end{theo}


We end this section with the following theorems that will be useful in the next section.
\begin{theo}{\cite[Theorem $2.17$]{ben1}} 
\label{lmsdlmsddmlslmsdlmsdl} 
For any  differentiable map $F:I\rightarrow {\sS'}(\R)$,  the element $<\hspace{-0.2cm}.,F(t)\hspace{-0.2cm}>$ is a differentiable stochastic distribution process which satisfies the equality:  
\vspace{-2ex}
\begin{equation*}
 \frac{d}{dt} <.,F(t)> \ = \  <.,\frac{d F}{dt} (t)>
\end{equation*}
\end{theo}

\begin{theo}{\cite[Lemmas $1$ and $2$ p.$73$-$74$]{GeCh2}}
\label{peogi}
Let $I$ be  an interval of $\R$, $t \mapsto F(t)$ be a map  from $I$ into ${\sS}_{-p}(\R)$,  $t \mapsto \varphi(t,.)$ be a map from $I$ into $\sS(\R)$ and $t_0 \in I$. If both maps $t \mapsto F(t)$ and  $t \mapsto \varphi(t,.)$ are continuous (respectively differentiable) at $t_0$,  then the function $t \mapsto <F(t),\varphi(t,.)>$ is continuous (respectively differentiable) at $t_0$. In this latter case  we have 
\begin{equation*}
 \frac{d <F(t),\varphi(t,.)> }{d t}(t_{0})  =  < \frac{dF}{dt}(t_{0}), \varphi(t,.)>  + < F(t),\frac{d}{dt}[ t \mapsto \varphi(t,.)](t_{0})>.
\end{equation*}
%
\end{theo}


%
\subsection{Stochastic integral with respect to Gaussian process}
\label{sjjjj}
Denote $G:=(G_t)_{t\in \rR}$ the process defined, for every $t$ in $\rR$, 
by $G_{t} :=  <.,g_{t}>$, where ${(g_t)}_{t\in\rR}$ is a family of functions of $L^2(\R)$. Denote $(t,s)\mapsto R_{t,s}$ the covariance function of $G$. We hence have 
$R_{t,s}:=\mathbf{E}[G_t\ G_s]={<g_{t},g_{s}>}_{L^{2}(\R)}$, for every $
(s,t)$ in $\rR^{2}$. We will note in the sequel $R_{t}$ instead of $R_{t,t}
$. 
$G$ is a Gaussian process which fulfills 
For the sake of notational simplicity we
can and will assume that $G_0 
\stackrel{{\textit a.s.}}{=}0$. Moreover, when the Gaussian process $G$ 
will admit a continuous modification, we will systematically use it and still 
call it $G$.
%
\subsubsection*{White Noise derivative of $G$}
Define the map  $g:\rR\rightarrow {\sS'}(\R)$ by setting 
$g(t) := g_t$. When $g$ is differentiable at point $t$, one 
denotes $g '_t$ its derivative. Denote $\lambda$ the Lebesgue measure on $\rR$ and define $L^1_{l\widetilde{oc}}
(\rR):=\{f:\rR\rightarrow \R \text{ is measurable }; \ f\in L^1((a,b)), \hspace{0.1cm} \text{for all finite interval  } (a,b) \textit{ s.t. } [a,b] \subset \rR \}$. 
{In this section and in the next one, we will make the following assumption:
\vspace{-4ex}

\makeatletter\tagsleft@true\makeatother 
\vspace{4ex}
\vspace{-6ex}

\textcolor{white}{
\begin{subequations}
 \label{A}
 \begin{align}
  a & = b \tag*{\text{($\mathscr{A}$)}} \label{Aa}\\
  a & = b \tag*{$\mathscr{A}_{\text{(i)}}$} \label{Ai}\\
  a & = b \tag*{$\mathscr{A}_{\text{(ii)}}$} \label{Aii}\\
  a & = b \tag*{$\mathscr{A}_{\text{(iii)}}$}\label{Aiii} \\
  a & = b \tag*{$\mathscr{A}_{\text{(iv)}}$} \label{Aiv}
 \end{align}
\end{subequations}
}

\vspace{-25ex}
\begin{align}
\label{Aba}
\tag*{(\hspace{-0.1ex}$\mathscr{A}_{\textcolor{white}{\text{iiiiii}\hspace{-1.6ex})}}$\hspace{-1.7ex})} 
\hspace{4ex} \begin{cases} 
(\text{i})  \text{ The map $g$ is continuous on } \rR,\\
(\text{ii})  \text{  The map $g$ is differentiable } \lambda-\text{almost everywhere on } \rR,  \\
(\text{iii}) \text{   There exists }  q   \text{ in } \N^*\text{ such that } {t\mapsto {|g'_t|}_{-q} \text{ belongs to } L^1_{l\widetilde{oc}}(\rR),}\\
(\text{iv}) \text{   For every } (a,b) \text{ in } {\rR}^{2} \text{   such that } a \leq b, \text{ one has, in $\sS'(\R)$, the equality: }\\
\textcolor{white}{pfojerpfoerofjerof} 
\end{cases} 
\end{align}

\vspace{-5.5ex}
\begin{equation}
\label{zodicjzoscjoifjs}
\tag{$E_{a,b}$}
g_{b} - g_{a} = \int^{b}_{a} \ g'_{u} \ du.
\end{equation}

%
%
%

\makeatletter\tagsleft@false\makeatother 
\vspace{1ex}

Define the set $\rR_{D}$ by setting $\rR_{D}:=\{t\in\rR; \vspace{1ex} \ g \text{  is differentiable at point } t\}$. 
For the sake of notational simplicity we will write  \ref{Ai}, 
(\text{resp.} \ref{Aii}, \ref{Aiii} or \ref{Aiv}), in the sequel, 
when one wants to refer to (i) (\text{resp.} to (ii), (iii) or (iv)) of Assumption\footnote{
An easy way to see if Assumption \Aa holds is to check if whether the following condition holds.
\vspace{-2ex}
$$
\hspace{-3ex}(\mathscr{D})
\begin{cases} 
(\text{i})  \text{ \small The map $g$ is continuous on $\rR$ and differentiable on  every finite interval } (a,b) \textit{ s.t. } [a,b] \subset \rR,\\
(\text{ii}) \text{  \small  There exists }  q   \text{ in } \N^*\textit{ s.t. } {t\mapsto {|g'_t|}_{-q} \text{ \small belongs to } L^1_{l\widetilde{oc}}(\rR).} 
\end{cases} 
$$
} \Aap.
\begin{rem}
\label{ozeufhioeruvheirezzezeedzedzedzdu}
 A first consequence of Assumption \ref{Aa} is that $g$ is 
``weakly'' locally absolutely continuous on $\rR$; this means that the map $t\mapsto <g_{t},\eta>$ is absolutely continuous on every finite interval $[a,b]$ of $
\rR$, for every $\eta$ in $\sS(\R)$.
\end{rem}

\begin{theo}[{\bfseries Gaussian White Noise \cite[Theorem-Definition 3.1]{JL17a}}]
\label{oijoifjsoidjdoijfsoijfsoijdsoifjdoisdjosdjf}
Define $W^{(G)}:=({W^{(G)}_{t})}_{t \in \rR_{D}}$ by setting: ${W^{(G)}}_{t} := \ <.,g'_t>$,
where the equality holds in  ${(\cS)}^*$. Then $W^{(G)}$ is 
is the ${(\cS)}^*$-derivative
of the process ${(G_t)}_{t \in \rR_{D}}$. We will sometimes note $\frac{dG_t}{dt}$ instead of $W^{(G)}_{t}$. Moreover  the map $t \mapsto {\| W^{(G)}_{t}\|}_{-p}$ is continuous if and only if $t \mapsto {|g'_t|}_{-p}$ is continuous. 
\end{theo}

\subsubsection*{Generalized functionals of $G$}
 In order to get 
a Tanaka formula in Section 
\ref{Tanaka}, 
we define here generalized 
functionals of $G$, by using \cite[Section $7.1$]{Kuo2}. We 
identify, here and in the sequel, any function $f$ of $L^{1}
_{loc}(\R)$ with its associated tempered distribution, denoted $T_f$, 
when it exists. In particular, one notes 
in this case: $<f,\phi>  \ = \int_{\R} f(t) \  \phi(t) \  dt$,  for 
every $\phi$ in  $\sS(\R)$. In this latter case we say that 
the tempered distribution $T:=T_{f}$ is of function type. 
Define the sets $\Z_R:=\{t \in \rR; \ R_{t} = 0\}$ and $
\Z^c_R:=\{t \in \rR; \ R_{t} > 0\}$.


%

%

\begin{theodef}
\label{lmqsdlmflkjdg}
 Let $F$ be a tempered distribution. For every $t$ in $\Z^c_R$, define 
 \vspace{-1.5ex}
\begin{equation*}
 F(G_t):= \frac{1}{\sqrt{2\pi {R}_{t}}} \sum\limits^{+\infty}_{k=0} \frac{1}{\ k! \ 
{{R}^{k}_{t}}}  <F, {\xi}_{t,k}> I_k\left(   {g_t }^{\otimes k} \right),
\end{equation*}
 where, for every $(x,k)$ in $\R\times \N$, ${\xi}_{t,k}(x) := {\pi}^{1/4} {(k!)}
^{1/2}  {{R}^{k/2}_{t}}\  {\exp}{\{-\frac{x^2}{4 {{R}_{t}}}\}}  {e}_{k}{(x/(\sqrt{2 
{R}_{t}}))}$. Then for all real $t$, $F(G_t)$ is a Hida distribution, called {\it 
generalized functional of $G_t$}.
 \end{theodef}

As stated in \cite{ben1}, when $F=f$ is of function type, $F(G_t)$ coincides with $f(G_t)$. 



%

%

\subsubsection*{S-Transform of $G$ and $W^{(G)}$}
The following theorem makes explicit the $S$-transforms of $G$, of the Gaussian white noise $W^{(G)}$ and of generalized functionals of $G$. Denote $\gamma$ the heat kernel density on $\R_+\times\R$ \textit{i.e.} 
\begin{equation}
\label{fctiongammadzedezdz}
\gamma(t,x):= \tfrac{1}{\sqrt{2\pi t}}\exp{\{ \tfrac{-x^2}{2t}  \}} \text{ if } t \neq 0 \text{ and } 0 \text{ if } t = 0. 
\end{equation}
\begin{theo}{\cite[Theorem 3.6]{JL17a}}
\label{tardileomalet}
For every $\eta$ in $\sS(\R)$ one has the following 
equalities:
\vspace{-1ex}
\bit
\iti  $S(G_t)(\eta) = {<g_t,\eta>}_{L^2(\R)}$, for every $t$ in $\rR$,
\vspace{-0.2cm}
\itii $S(W^{(G)}_{t})(\eta) = {<g'_t,\eta>} =\frac{d}{dt} [{<g_t,\eta>}_{L^2(\R)}] $, for every $t$ in  $\rR_{D}$;
\itiii  For $p \in \N$, $F \in \sS_{\hspace{-0.15cm}-p}(\R)$, and $t$ in $\Z^c_R$,  \  $S(F(G_t))(\eta) = \left< F, \gamma \left(R_{t}, .- {<g_t, \eta>}  \right) \right>$.
\smallskip

Furthermore, there exists a constant $D_p$, independent of $F,t$ and $\eta$, such that:
$$\forall t \in \Z^c_R, \hspace{1cm} \ {|S(F(G_t))(\eta)|}^2  \leq \ \max\{   {{R}^{-2p}_{t}}; \ {{R}^{2p}_{t}} \} \    {{R}^{-1/2}_{t}}   \  \ D_p  \ {|F|}^2_{-p} \  \exp\{{| \eta |}^2_{p}\}.$$

\eit
%
%
%
%
%
%
%
%
\end{theo}

\subsubsection*{The Wick-Itô integral with respect to Gaussian processes}
%
We still assume in this section that Assumption \Aa holds and still denote $I$ a Borel set of $\rR$. 
%
%
%
\begin{defi}[Wick-Itô integral \textit{wrt} Gaussian process]
\label{oezifhherioiheroiuh}
Let $X:\rR \rightarrow {(\cS)}^*$ be a process such that the process 
$t \mapsto X_t \diamond W^{(G)}_t$ is  \oS-integrable on $\rR$. 
The process $X$ is then said to be $dG$-integrable on $\rR$ (or 
integrable  on $\rR$), \textit{wrt} the Gaussian process $G$. The 
Wick-Itô integral of $X$ \textit{wrt} $G$, on $\rR$, is defined by 
setting:
\vspace{-1ex}

\begin{equation}\label{eigfrretth}
\int_{\rR} X_s \ \di G_s :=  \int_{\rR} X_s \diamond W^{(G)}_s \ ds.
\end{equation}

For any $I$ in $\cB(\rR)$, define $\int_{I} X_s \ \di G_s:= \int_{\rR} \ {\i1}
_{I}(s) \ X_s \ \di G_s$.
\end{defi}

 The Wick-Itô integral of an $ ({\cS}^*)$-valued process, \textit{wrt} $G$ is then an element of \oS. It is easy to see that Wick-Itô integration \textit{wrt} $G$, is linear and, assuming it belongs to $(L^2)$, centered. 

\begin{theo}{\cite[Theorem 3.13]{JL17a}}
\label{dlsdqsdhuqsgdudguqsgdqudgquygdsudgqsuydgqsuydgqsuyqdgq}
Let $p \in \N$, $[a,b]$ be an interval of $\Z^c_R$ and let $F:[a,b] \rightarrow 
\sS_{\hspace{-0.15cm}-p}(\R)$ be a continuous map. If $t\mapsto \max \{R^{-
p-1/4}_{t}; R^{p-1/4}_{t} \}$ belongs to  $L^1([a,b])$ (\textit{resp.} there 
exists an integer $q$  such that $t\mapsto {|g'_t|}_{-q}\max \{R^{-p-1/4}
_{t}; R^{p-1/4}_{t} \}$ belongs to  $L^1([a,b])$), then the stochastic 
distribution process $F(t,G_t)$ is ${(\cS)^*}$-integrable  (\textit{resp.} $dG$-integrable) on $[a,b]$ (\textit{resp.} on $(a,b)$).
\end{theo}

\vspace{-2ex}
\subsection{Gaussian Processes in $\mathscr{G}$ of ``reference''}
\vspace{-1ex}
To see the generality of the results on local times we present here, we will consider in this paper, classical Gaussian processes, made with elements of $\mathscr{G}$. These processes are: Brownian motion and  
Brownian bridge, fractional and multifractional Brownian motions as 
well as $\mathscr{V}_{\gamma}$ - processes.
\smallskip

\textit{\bfseries Fractional, Multifractional Brownian motions and $\mathscr{V}_{\gamma}$ - processes}
\smallskip

Readers interested in an exhaustive presentation of fBm or mBm may refer to \cite{Nu} for fBm and to \cite{LLVH} for mBm, as well as to the references therein. Recall that fBm, which was Introduced in \cite{Kol} and popularized in \cite{Mandelbrot1968}, is a centered Gaussian process, the covariance function of which is denoted $R_{H}$ and is given by: 
\vspace{-1ex}
\begin{equation*}
  R_{H}(t,s) :=  \frac{1}{2} ( {|t|}^{2{H}} + {|s|}^{2{H}} - {|t-s|}^{2{H}}),
\end{equation*}
where $H$ belongs to 
$(0,1)$, and is usually called the Hurst exponent. When $H=1/2$, fBm reduces to standard 
Brownian motion. 
MBm, which is a Gaussian extension of fBm, was 
introduced in \cite{PL} and in \cite{ABSJDR} in order to 
match any prescribed non-constant deterministic local 
regularity and to decouple this property from long range 
dependence. 
A mBm 
on $\R$, with functional parameter $h:\R \rightarrow (0,1)$, 
is a Gaussian process $B^h:={(B^h_t)}_{t\in \R}$ defined, 
for all real $t$, by $B^h_t:=\sifbm(t,h(t))$, where  $
\sifbm:={(\sifbm  (t,H))}_{(t,H) \in \R\times(0,1)}$ is 
fractional Brownian field on $\R\times(0,1)$ (which means 
that  $\sifbm$ is a Gaussian field, such that, for every $H$ 
in $(0,1)$, the process ${(\sifbm  (t,H))}_{t \in \R}$ is a fBm 
with Hurst parameter $H$). In other words, a mBm is 
simply a ``path'' traced on a fractional Brownian field. 
Note also that when $h$ is constant, mBm reduces to fBm. 
For any deterministic function $\gamma:\R_{+}\rightarrow \R$, $
\mathscr{V}_{\gamma}$ - processes 
are defined in \cite[Prop. 1]{MV05} as being the processes, denoted ${\widetilde{B}}^{\gamma}:= {({\widetilde{B}}^{\gamma}_{t})}_{t\in [0,T]}$, by setting: 
\vspace{-2ex}
\begin{equation}
\label{erofkorekdeddeddeezzeep}
{\widetilde{B}}^{\gamma}_{t}:=\int^{t}_{0} \varepsilon(t-s) \ dW_{s}; \hspace{0.5cm} \forall t\in [0,T],
\end{equation}
with $\gamma:\R_{+}\rightarrow \R$ such that $
\gamma^{2}$ is of class $C^{2}$ everywhere in $\R_{+}$ 
except in $0$; and such that ${(\gamma^{2}})'$ is non increasing. The map $\varepsilon:\R^{*}_{+}\rightarrow \R$ is 
defined by setting $\varepsilon:= \sqrt{{(\gamma^{2}})'}$. Subset of $\mathscr{G}$, the set 
$\mathscr{V}_{\gamma}$ contains Gaussian processes, that can be 
more irregular than any fBm. However it does not contain fBm (nor mBm) since $\mathscr{V}_{\gamma}$ only contains processes the regularity of which
 remains constant along the time).
\smallskip

\textit{\bfseries Operators ${(M_{H})}_{H\in(0,1)}$}
\smallskip

The operator $M_{H}$ will be useful in the sequel, not only to provide one with a representation of fBm and of mBm under the form \eqref{erofkorekp}, but also to verify that Assumption \Aa made page \pageref{Aa}, hold for both fBm and mBm. 
Let $H$ belongs to $(0,1)$; following \cite{ell} and \cite[Section $2.2$]{JLJLV1}, define the $L^{2}(\R)$-valued operator $M_H$, in the Fourier domain by: $\widehat{M_H(u)}(y) := \tfrac{\sqrt{2\pi}}{c_H} \hspace{0.1cm} |y|^{1/2-H} \ \widehat{u}(y), \hspace{0.5cm} \forall y \in \R^*$,
where $c_{x}$ is defined, for every $x$ in $(0,1)$ by $c_x:=  {\big(\tfrac{2\pi}{\Gamma(2x+1) \sin(\pi x)} \big)}^{\frac{1}{2}}$.
This operator is well defined on the homogeneous Sobolev space
$ L^2_H(\R):= \{u \in {\sS'}(\R) \hspace{0.1cm}:\hspace{0.1cm} \widehat{u} = T_f; 
\hspace{0.1cm} f \in L^{1}_{loc}(\R) \hspace{0.1cm}\textnormal{and}
\hspace{0.1cm} \| u \|_H < +\infty \},$
%
%
%
where the norm $\ds \|\cdot\|_H$ derives from the inner product denoted ${\langle \cdot , \cdot \rangle}_H$, which is defined on $\ds L^2_H(\R)$ by: ${\langle u,v \rangle}_H := \frac{1}{c^2_H} \int_\R |\xi|^{1-2H} 
{\widehat{u \hspace{0.1cm}}(\xi) }\ \overline{{\widehat{v \ }
(\xi) }} \ d{\xi}$.
$M_{H}$ being an isometry from $(L^2_H(\R), \ds \|\cdot\|_H)$ 
into $(L^2(\R), \ds \|\cdot\|_{L^2(\R)})$, it is clear that, for every $(H,t,s)$ in $(0,1)\times \R^{2}$, ${<M_{H}(\i1_{[0,t]}), M_{H}(\i1_{[0,s]})>}_{L^{2}(\R)}= R_{H}(t,s)$.
%
%
%
We will say that an mBm is normalized when its covariance function, denoted $R_{h}$, verifies the equality:
\begin{equation}
\label{azazazaz}
R_{h}(t,s)=     \tfrac{ c^2_{h_{t,s}}}{c_{(h(t))}c_{(h(s))}}  \ \big[\tfrac{1}{2} \big( {|t|}^{2h_{t,s}} + {|s|}^{2h_{t,s}}  - {|t-s|}^{2h_{t,s}} \big)\big],
\end{equation}
where $h_{t,s} := \frac{h(t)+h(s)}{2}$ and $c_{x}$ has been above, right after $\widehat{M_H(u)}(y)$.

\begin{ex}[Gaussian Processes in $\mathscr{G}$ of ``reference'']
\label{oevieroivjedzezdzeddo}
Let $H$ be real in $(0,1)$ and $h:\R \rightarrow (0,1)$ be a deterministic measurable function. Define the processes 
\vspace{-0.5cm}

\begin{align*}
&B:=\{   <.,\i1_{[0,t]}>;\ t\in\R   \};& \
&\widehat{B}:=\{<.,\i1_{[0,t]} - t\cdot \i1_{[0,1]}>;\ t\in[0,1]   \};&\\
&B^{H}:=\{ <.,M_{H}(\i1_{[0,t]})>;\ t\in\R   \};& \
& B^{h}:=\{ <.,M_{h(t)}(\i1_{[0,t]})>;\ t\in\R   \};\ & \\
&{\widetilde{B}}^{\gamma}:=\{ <.,\i1_{[0,t)} \cdot \varepsilon(t-.))> \text{ if }  t\in\R^{*}_{+}  \ \& \  {\widetilde{B}}^{\gamma}_{0}:= 0 \}.& \
\end{align*}
%
We know, thanks to Section \ref{cddcwcdsd}, that $B$ is a 
Brownian motion on $\R$, that $\widehat{B}$ is Brownian bridge on $[0,1]$, 
that  ${B}^{H}$ is a fBm of Hurst index $H$, 
that  $B^{h}$ is a normalized mBm of functional parameter 
$h$ and 
that ${\widetilde{B}}^{\gamma}$ is a $\mathscr{V}_{\gamma}$ - process.
\end{ex}
A word on notation: $B^H_.$ or $B^{h(t)}_.$ will always 
denote an fBm with Hurst index $H$ or $h(t)$, while 
$B^h_.$ will stand for an mBm. 

\section{Itô Formulas for generalized functionals of $G$}
\label{Itô}
 The main result of this section is an Itô Formula, in $(L^{2})$-sense, for 
generalized functionals of $G$ (Theorem 
\ref{zpzpkpskpsokpdokposdkfpsdokfpsokshfigushfzdiziyuezaiozo}). This latter  
will not only provide us with a Tanaka formula, 
but will also give us a precise reason to define the weighted local time of 
$G$, the way we do it, in Section \ref{Tanaka}. To do so
we first need to establish an Itô formula in \oS. This is done in Subsection 
\ref{zpepoko}.


{
Let $[a,b]$ be an interval of $\R$ and $j:[a,b]
\rightarrow \R$ be a function of bounded variation. Denote $\alpha_{j}$ the signed measure such that $j(t)=\alpha_{j}([a,t])$, for every $t$ in $[a,b]$. For any function $f:[a,b]\rightarrow \R$,  denote $ \int^{b}_{a} f(s) \ dj(s)$ or $ \int^{b}_{a} f(s) \ d\alpha_{j}(s)$ the \textit{Lebesgue-Stieljes integral} of $f$ with respect to $j$, assuming it exists. In this latter case, we will write  that $f\in L^{1}(I,dj(t))$ or $L^{1}(I,\alpha_{j})$.}
%
%
%
In the particular case where the function $f$ is continuous on $[a,b]$, the \textit{Lebesgue-Stieljes integral} of $f$ exists and is also equal to the \textit{Riemann-Stieljes integral} of $f$, which is denoted and defined by:
\begin{equation}
\label{pzodkk}
\text{(R.S.) }  \int^{b}_{a}\ f(s)\  d j(s):= \underset{\pi \to 0}{\lim}  \sum^{n}_{i=1}\ f(\xi^{(n)}_{i}) \ (j(x_{i})-j(x_{i-1})),
\end{equation}
where the convergence holds uniformly on all finite partitions ${\mathscr {P}}^{(n)}_{\pi}:=\{a:=x_{0} \leq x_{1} \leq \cdots \leq x_{n}:=b\}$ of $[a,b]$ such that $\underset{1\leq i \leq n}{\max} (x_{i}-x_{i-1}) \leq \pi$ and such that $\xi^{(n)}_{i}$ belongs to $[x_{i-1},x_{i}]$. 
The following result, will be used extensively in the sequel of this section.

\begin{lem}{\cite[Lemma 4.1]{JL17a}}
\label{zpodkzpeodkzep}
Let $[a,b]$ be a finite interval of $\R$, $I$ (\textit{resp.} $J$) an interval of  $\R_{+}$ (\textit{resp.} of $\R$) and let $L:[a,b]\times I\times J$ be a $C^{1}$-function. Let $f:[a,b]\rightarrow I$ and $j:[a,b]\rightarrow J$ be two continuous functions of bounded variation on $[a,b]$. Then one has the following equality:
\vspace{-3ex}
\begin{align}
 \label{zpdozepkddzzdzdzdpk}
 L(b,f(b),j(b)) -  L(a,f(a),j(a)) &= \int^{b}_{a}\    \frac{\partial L}{\partial u_{1}}(s,f(s),j(s))  \ ds +  \int^{b}_{a}\    \frac{\partial L}{\partial u_{2}}(s,f(s),j(s))  \ df(s) \notag \\
&+ \int^{b}_{a}\    \frac{\partial L}{\partial u_{3}}(s,f(s),j(s))  \ dj(s).
\end{align}
 
\end{lem}
%
%


In the remaining of this paper, and unless otherwise specify,
the measure $m$ denotes a measure, that may be $\sigma$-finite or signed.

\begin{lem}
\label{podiuhyiuuhiudokdfzdokfdepsoksdpokfsdpokfspokdpfoksdpfskfodkof}
Let $T>0$ and $v:[0,T]\times \R \rightarrow \R$ be a continuous function 
such that there exists  a couple $(C_T,\lambda_T)$ of
$\R\times\R^*_+$ such that $\underset{t \in [0,T] }{\max}{|v(t,y)|} \leq C_T \
e^{\lambda_T y^2}$ for all real $y$. Define; for every $a>\lambda_T$, the map  $J_v:\R_+ \times (0,1/4a) \times \R \rightarrow \R$ by 
setting: 
\begin{equation}
\label{osicosicjsdoicj}
J_v(t,u_{1},u_{2}):= \int_{\R} \  v(t,x)\cdot \gamma(u_{1},x-u_{2}) \  dx.  
\end{equation}
Then $J_v$ is 
well defined. Moreover \hspace{-0.5cm} $\lim\limits_{(t,u_{1},u_{2}) \to (t_{0},0^+,l_0)}  
\hspace{-0.15cm} J_{v}(t,u_{1},u_{2}) = v(t_{0},l_{0})$, $\forall\  (t_{0},l_0)$ in $[0,T]\times \R$.
\end{lem}

\subsection{Itô Formula  in \oS \ for generalized functionals of $G$ on an interval of $\Z^c_R$}
\label{zpepoko}

 Before establishing an Itô formula in \oS, let us first fix some notations.
%
%
For a tempered distribution $G$ and a positive integer $n$, let $G^{(n)}$ denote the $n^{\text{th}}$ distributional derivative of $G$. We also write $G':= G^{(1)}$. Hence, by definition, the equality  $<G',\varphi>=-<G,\varphi'>$ holds for all $ \varphi$ in $\sS(\R)$.
For a map $t \mapsto F(t)$ from $[a,b]$ to $\sS_{\hspace{-0.15cm}-p}(\R)$ we will note  $\frac{\partial^n F }{\partial x^n}(t)$ the quantity ${(F(t))}^{(n)}$, that is the $n^{\text{th}}$ derivative in $\sS^{'}(\R)$, of the tempered distribution $F(t)$.
Hence we may consider the map  $t \mapsto \frac{\partial^n F }{\partial x^n}(t)$ from $[a,b]$ to $\sS'(\R)$. Moreover for any $t_0$ in $[a,b]$, we will note $\frac{\partial F}{\partial t}(t_0)$ the quantity $\lim\limits_{r \to 0} \frac{F(t_{0}+r)-F(t_0)}{r}$, when it exists in $\sS_{\hspace{-0.15cm}-p}(\R)$, for a certain integer $p$. When it exists, $\frac{\partial F}{\partial t}(t_0)$ is a tempered distribution, which is said to be the derivative of the distribution $F(t)$, with respect to $t$ at point $t=t_0$.
In line with Subsection \ref{sjjjj}, we then  define, for $t_0$ in $[a,b]$ and a positive integer $n$,  the following quantities: $\frac{\partial^n F }{\partial x^n}(t_0,G_{t_0}) := {(F(t_0))}^{(n)} (G_{t_0})$ and $\frac{\partial F }{\partial t}(t_0,G_{t_0}) := \left(\frac{\partial F }{\partial t}(t_0)\right) (G_{t_0})$.

\smallskip

%
%
%
%

%
%

\begin{theo}
\label{ziejoezijoez}
Let $[a,b]$ be an interval of $\Z^c_R$, $p \in \N$, and let $F$ be an element of $\cC^1([a,b],\sS_{\hspace{-0.15cm}-p}(\R))$ such that both maps $\frac{\partial F}{\partial x}$ and $\frac{\partial^2 F}{\partial x^2}$, from $[a,b]$ into $\sS_{\hspace{-0.15cm}-p}(\R)$, are continuous.
Assume moreover that \Aa holds and that the map $t\mapsto R_{t}$ is both continuous and of bounded variations on $[a,b]$. Then the following equality holds in \oS:

\vspace{-3ex}

\begin{equation}
\label{orijvedze2}
F(b,G_{b}) - F(a,G_{a}) = \int^b_a \ \tfrac{\partial F}{\partial t}(s,G_{s}) \ ds + \int^b_a \ \tfrac{\partial F}{\partial x}(s,G_{s}) \ \di G_{s}
+ \tfrac{1}{2} \ \int^b_a \ \tfrac{\partial^2 F}{\partial x^2}(s,G_{s}) \  dR_{s}.
\end{equation}
\end{theo}

%


\begin{pr}
Assumption \Aaiii and the fact the map $t\mapsto R_{t}$ is upper and 
lower bounded on $[a,b]$ allow one to apply Theorem 
\ref{dlsdqsdhuqsgdudguqsgdqudgquygdsudgqsuydgqsuydgqsuyqdgq}  to $
\tfrac{\partial F}{\partial t}(s,G_{s})$, $\tfrac{\partial^{2} F}{\partial x^{2}}
(s,G_{s})$ (with the Borel measure defined by $t\mapsto R_t$),  and to $\tfrac{\partial F}{\partial x}(s,G_{s})$. 
This and \cite[Remark 6]{JL17a} 
entail that all the integrals on the right side of 
\eqref{orijvedze2} exist.
Thanks to Lemma \ref{dkdskcsdckksdksdmksdmlkskdm} it is then sufficient to show the equality of the $S$-transforms of both sides of \eqref{orijvedze2}. Using Theorem \ref{tardileomalet} one then just has to establish, for every $\eta$  in $\sS(\R)$, the following equality:
\vspace{-0.5cm}

\begin{align}
\label{pzodkpzdkpzeodkpz}
I_{a,b}(\eta) &= \hspace{-0.1cm}\int^{b}_{a}\hspace{-0.2cm} <\tfrac{\partial F }{\partial t}(s),\gamma (R_{s}, .- {<g_s, \eta>}) >  ds +   \int^{b}_{a}\hspace{-0.2cm}  < \tfrac{\partial F}{\partial x }(s), \gamma(R_{s}, .- {<g_s, \eta>}) > <g'_{s},\eta>   ds \notag\\														& \hspace{0.15cm}  +   \int^{b}_{a}  \ <  \tfrac{\partial^{2} F}{\partial x^{2} }(s), \gamma(R_{s}, .- {<g_s, \eta>}) > \ dR_{s},
\end{align}

where  we set $I_{a,b}(\eta):=<F(b),\gamma(R_{b}, .- {<g_b, \eta>})>-<F(a),\gamma(R_{a}, .- {<g_a, \eta>})>$. Define the map $L:[a,b]\times \R_{+}\times\R$ by setting: $L(u_{1},u_{2},u_{3}):=<F(u_{1}),\gamma(u_{2},.-u_{3})>$.
The assumptions made on $F$ entail that $L$ is a $C^{1}$-function. Moreover, using the equality $\tfrac{\partial \gamma}{\partial t} =  \tfrac{1}{2} \tfrac{\partial^2 \gamma}{\partial x^2}$, one gets:
\begin{align*}
 < \tfrac{\partial F}{\partial x }(s), \gamma(R_{s}, .- {<g_s, \eta>}) >   &=  \ < F(s), \tfrac{\partial \gamma}{\partial x }(R_{s}, .- {<g_s, \eta>}) >, \notag \\
 <  \tfrac{\partial^{2} F}{\partial x^{2} }(s), \gamma(R_{s}, .- {<g_s, \eta>}) > &= \tfrac{1}{2} < F(s), \tfrac{\partial \gamma}{\partial t }(R_{s}, .- {<g_s, \eta>}) >.
\end{align*}
%
%

Since  $I_{a,b}(\eta) = L(b,R_{b},\gamma(<g_{b},\eta>)) - L(a,R_{a},\gamma(<g_{a},\eta>))$,  \eqref{pzodkpzdkpzeodkpz}  also reads as \eqref{zpdozepkddzzdzdzdpk}
with $f(s):=R_{s}$ and $j(s):=j_{\eta}(s):=<g_{s},\eta>$. According to Remark \ref{ozeufhioeruvheirezzezeedzedzedzdu}, $j_{\eta}$ is absolutely continuous on $[0,T]$. Lemma \ref{zpodkzpeodkzep} therefore applies and 
achieves the proof.
\end{pr}
%
%
%
%
%
%
%

%

\begin{rem}
As Theorem 
\ref{zpzpkpskpsokpdokposdkfpsdokfpsokshfigushfzdiziyuezaiozo} below will 
show, one can extend \eqref{orijvedze2} to the case where $[a,b]\cap \zr 
\neq \emptyset$, if one makes additional assumptions on $R$ and $F$. However, we will not extend \eqref{orijvedze2} directly since one does not need 
such a general formula to get Theorem 
\ref{zpzpkpskpsokpdokposdkfpsdokfpsokshfigushfzdiziyuezaiozo}, 
and hence a Tanaka formula. 
\end{rem}


\subsection{Itô Formula in $(L^{2})$ for certain generalized functionals of $G$ on any interval of $\rR$}
\label{osozirjorjoeric}
To get a Tanaka formula, 
one needs both: to generalize \eqref{orijvedze2} to the case where $[a,b]$ 
contains elements of $\Z_R$, and to specify \eqref{orijvedze2} to the case 
where $F(t,x) =F(x)$. However, Theorem \ref{ziejoezijoez} does not extend immediately to the case where $[a,b]$ intersects the set $\Z_R$ because the generalized functionals of G 
are not defined in this situation. As in \cite{ben1}, we now extend our previous results to deal with this difficulty. 
%
Having in mind the definition of distributions of function type, given right above Theorem-Definition \ref{lmqsdlmflkjdg}, one has the following result.

%
%
%
%

\begin{theo}
\label{zpzpkpskpsokpdokposdkfpsdokfpsokshfigushfzdiziyuezaiozo} 
Let $[a,b]$ be an interval of $\rR$, $F$ be a distribution of function type that is 
continuous on $\za:=[a,b]\cap\Z_R$.  Assume that the first distributional derivative 
of $F$, noted $F'$, is also of function type and that $F$ and $F'$ both have 
polynomial growth. Assume moreover that \Aa holds, that the map $t\mapsto R_{t}$ is both continuous and of bounded variations on $[a,b]$ and that the following assumptions are satisfied:
\vspace{-0.25cm}
 
\bit
\iti$\lambda(\za) = \alpha_{R}(\za)=0$,
\itii $t\mapsto R^{{\scriptscriptstyle-1/2}}_{t}$ belongs to $L^1([a,b], dR_{t})$,
\itiii $\exists \ q \in \N$ such that $t\mapsto {|g'_t|}_{-q} \ R^{-1/2}_{t} \in L^1([a,b], dt)$,
\eit
 
where $\lambda$ denote the Lebesgue measure on $\rR$ and where the signed measure $\alpha_{R}$ has been defined at the beginning of Section \ref{Itô}. Then, the 
following equality holds in $(L^2)$:

%

\vspace{-0.15cm}
\begin{equation}
\label{ozoisdjcjodsjcdsfjciojciosd}
F(G_{b})- F(G_a) = \int^b_a F'(G_{t}) \ \di G_{t} + \frac{1}{2} \ \int^b_a \ F''(G_{t}) \ dR_{t}.
\end{equation}
\end{theo}

\begin{pr}
Let us first show that $t\mapsto F'(G_{t}) \hspace{0.2ex}\diamond \hspace{0.2ex} W^{(G)}_t$ is \oS-
integrable on $[a,b]$ and that $t\mapsto F''(G_{t})$ is \oS-integrable on $
[a,b]$ with respect to $\alpha_{R}$. Denote $\wa:=\{t\in [a,b]; \ R_{t}>0\}$.
\vspace{-0.25cm}

\begin{itemize}
\ita \textit{$t\mapsto F'(G_{t}) \diamond W^{(G)}_t$ is \oS-integrable on $(a,b)$.}
\end{itemize}

Denote  $\overline{R}_{{\scriptstyle a,b}}:= \underset{t\in [a,b]}{\sup}{R_{t}}$ and $\widehat{D}:= \sup {\{\ {|{e}^2_{k}(u)|}: k \in \N, u \in \R \}}$.
We have the following upper bound, valid for all $k$ in $\N$ and $t$ in $\wa$,
\vspace{-0.5cm}

\begin{align*}
\label{oweifjwoeifjweoifweoifjdwewedwedwedwe}
{|<F', {\xi}_{t,k}>|}^2  = {|<F, {\xi}^{'}_{t,k}>|}^2 =  {\pi}^{1/2} {k!} \ R^{k}_{t} \ {|<F, \tfrac{d}{dx} [e_{k}(x/\sqrt{2} R^{1/2}_{t})  \exp\{-x^{2}/4 R_{t}\}]>|}^2,
\end{align*}

where ${\xi}^{'}_{t,k}$ denotes $\frac{d}{dx} [x\mapsto {\xi}_{t,k}(x)]$ (which has been defined right after Theorem-Definition \ref{lmqsdlmflkjdg}). The right hand side of the previous equality can be decomposed into two parts. Since we get the upper bound of these two parts similarly, we only show here how to get the upper bond of 
\vspace{-0.25cm}

\begin{equation*}
{|<F,  {(2\ R_{t})}^{{\scriptstyle -1/2}} e'_{k}(x/\sqrt{2} R^{1/2}_{t})  \exp\{-x^{2}/4 R_{t}\}>|}=: J_{t,k}. 
\end{equation*}

Denote $\widetilde{K}$ (\textit{resp.} $N$) the real  (\textit{resp.} the integer) such that $|F(x)|\leq \widetilde{K} (1+{|x|}^{N})$, for every real $x$. Using the relation $e'_{k}(x)=\sqrt{\frac{k}{2}} e_{k-1}(x)- \sqrt{\frac{k+1}{2}} e_{k+1}(x)$ (see \cite[p.$354$]{Kuo2}), valid for all positive integer $k$, and the polynomial growth of $F$, one immediately gets:

\vspace{-0.5cm}

\begin{align*}
J_{t,k} &\leq {(2\ R_{t})}^{{\scriptstyle -1/2}}  \int_{\R} |F(x)| \ \sqrt{2(k+1)} \ \widehat{D}\ e^{-x^{2}/4 R_{t}} \ dx
 \notag\\ 
&\leq     {(1+\sqrt{\overline{R}})}^{N} \widehat{D}  \widetilde{K}  \int_{\R} {(1+|u|)}^{N} \ 
e^{-u^{2}/4} \ du\ \sqrt{k+1}.
\end{align*}

It is now easy to obtain the next inquality, for a real $D'$ which does not depend on $k$ nor $t$:

\vspace{-0.25cm}

\begin{equation}
\label{oerifjoe}
\forall (t,k)\in \wa \times \N, \hspace{0.5cm} {|<F', {\xi}_{t,k}>|}^2   \leq D' (k+1)! \ R^{k}_{t}.\\
\end{equation}

Moreover, since the norm operator of $A^{-1}$ is equal to $1/2$ 
(see \cite[p.$17$]{Kuo2}), we get for a real $D$ that does not 
depend on $t$, using Theorem-Definition \ref{lmqsdlmflkjdg} and 
\eqref{oerifjoe},  for all $t$ in $\wa$:

\vspace{-0.15cm}

\begin{equation}
\label{oweifjwoeifjweoifweoifj}
{\|F'(G_{t})\|}^2_{-1} =  \tfrac{1}{2\pi R_{t}}  {\scriptstyle{\sum\limits^{+\infty}_{k=0}}} {(k!)}^{-1} \
 R^{-2k}_{t} \ {<F', {\xi}_{t,k}>}^2 \ {|A^{-1}g_t|}^{2k}_0 \leq D \cdot R^{-1}_{t}.
\end{equation}

In view of Lemma \ref{dede} and of \eqref{oweifjwoeifjweoifweoifj}, it is clear that there exists a real $K$, such that, for every $\eta$ in $\sS(\R)$ and $t$ in $\wa\cap \rR_{D}$, $| S(F'(G_{t}) \diamond W^{(G)}_t)(\eta)| \leq   {\|F'(G_{t})\|}_{-1} \ {\|W^{(G)}_t\|}_{-q} \  e^{{|\eta|}^2_{q} }  \leq  K  \ {|g'_t|}_{-q}\ R^{-1/2}_{t}\ e^{{|\eta|}^2_{q} }$.
One can then apply Thm. \ref{peodcpdsokcpodfckposkcdpqkoq} and conclude that $\int^b_{a} F'(G_{t}) \ \di G_{t}\in$ \oS.
\vspace{-1ex}

\begin{itemize}
\itb \textit{$t\mapsto F''(G_{t})$ is \oS-integrable on $[a,b]$, with respect to $\alpha_{R}$.}
\end{itemize}
\vspace{-0.35cm}

Since $F'$ also has polynomial growth we get, by replacing $F$ by $F'$ in  \eqref{oweifjwoeifjweoifweoifj}, using Lemma \ref{dede}, $| S(F''(G_{t}))(\eta)| \leq  K  \  R^{-1/2}_{t} \ e^{{|\eta|}^2_{1}}$, for every $t$ in $\wa$,
Thanks to Assumption $(i)$, we know that the right hand side of the previous inequality belongs to $ L^1([a,b]\backslash\za, dR_{t})$. Theorem 
\ref{peodcpdsokcpodfckposkcdpqkoq} then applies and ensures that $ \int^b_a \ F''(G_{s}) \ dR_{s}$ exists in \oS.
%
%
%
%
%
%
%
%
%
%
%
To establish \eqref{ozoisdjcjodsjcdsfjciojciosd}, it is then sufficient 
to show that the following equality holds in $(L^{2})$:
\vspace{-0.15cm}

\begin{equation}
\label{ozdzedezrijvpeofperdzzdzofkr2}
F(G_{b}) - F(G_{a}) = \int_{\wa}   F'(G_{t})\ \di G_{t} + \frac{1}{2} \ \int_{\wa}  F''(G_{t}) \ dR_{t}.
\end{equation}

Having in mind $(iii)$ of Theorem \ref{tardileomalet}, we proceed as in the proof of \cite[Theorem 4.4]{JL17a}. Since $\wa$ is an open set of $[a,b]$, it can be written under the form:
\vspace{-0.15cm}
 
\begin{equation}
\label{ocdedcdspcokdspockzpcozdzezdodjz}
\wa=[a,c)\ {\sqcup} \underset{i\in \N}{\bigsqcup} (c_{i}, d_{i})\ {\sqcup}\ (d,b],    
\end{equation}
 
where all the intervals in \eqref{ocdedcdspcokdspockzpcozdzezdodjz} are 
disjoint and where, by convention, $(x,y)=(x,y]=\emptyset$, for every reals 
$x$ and $y$ such that $x\geq y$. Note moreover that every element of $\{c_{i},d_{i},\ i\in \N \}$ (as well as $c$, if $[a,c)\neq\emptyset$, and $d$, if $
(d,b] \neq\emptyset$)  belongs to $\za$. Exactly as in the proof of \cite[Theorem 4.4]{JL17a}, one needs to distinguish between two cases:
 \smallskip

{\bfseries First case:} $\exists\ (c',d') \in {[a,b]}^{2}$ with $c'<d'$ \textit{s.t.} $[a,c')$ and $(d',b]$ are both subsets of $\wa$.
\vspace{0.1cm}

\smallskip

{\bfseries Second case:}  $\nexists\ (c',d')\in {[a,b]}^{2}$ with $c'<d'$ \textit{s.t.} $[a,c')$ and $(d',b]$ are both subsets of $\wa$.
\vspace{0.1cm}

 In both cases, the same reasoning as the one we used in the first case of the proof of \cite[Theorem 4.4]{JL17a}
  (with Theorem \ref{ziejoezijoez} instead of Lemma \ref{zpodkzpeodkzep} and Lemma \ref{podiuhyiuuhiudokdfzdokfdepsoksdpokfsdpokfspokdpfoksdpfskfodkof} with $v(t,.)=F$) 
%
%
%
%
%
%
%
%
applies and allows us to conclude.
\end{pr}

\begin{rem}
\label{ozeidjezodijzeodjezodjezodizeofnzeifubzeodazkpdoajodizb}
As it has been noticed in \cite[Remark 4.6]{ben1} (in the case of fBm), Equality \eqref{ozoisdjcjodsjcdsfjciojciosd} holds in $(L^2)$ in the sense that both sides belong to $(L^2)$. If $F''(G_{s})$ belongs to $(L^2)$ and under the additional assumption that $ \int^v_u \ \| F''(G_{s})\|_0 \ dR_{s} <+\infty$, then  $\int^v_u \ F''(G_{s}) \ dR_{s}$ is an $(L^2)$-valued Bochner integral (see Appendix \ref{appendiceB} for the definition of Bochner integral) and all members of both sides of\eqref{ozoisdjcjodsjcdsfjciojciosd} belong to $(L^2)$.
Note moreover that the assumption $t\mapsto R_{t}$ is bounded on $[a,b]$ is sufficient but not necessary; however, if one wants to relieve assumption on $R$, one will have to weigh down the ones on $F$, $F'$ or on the process $G$. 
\end{rem}
 Once again, all the Gaussian processes of ``reference'' fulfill assumptions of 
Theorem 
\ref{zpzpkpskpsokpdokposdkfpsdokfpsokshfigushfzdiziyuezaiozo}. The 
only non obvious case is the one of $V_{\gamma}$ - processes. In this latter case, 
$t\mapsto R_{t}$ is absolutely continuous on $[0,b]$ and therefore condition $(iii)$ of 
Theorem 
\ref{zpzpkpskpsokpdokposdkfpsdokfpsokshfigushfzdiziyuezaiozo} is fulfilled. 
Besides, because $R$ is: absolutely continuous on $[0,b]$, non 
decreasing and such that $\lim_{u\to 0+} R'_{u} = +\infty$, it is easy to  show, using \cite[(B.6)]{JL17a}, that
condition $(i)$ of 
Theorem 
\ref{zpzpkpskpsokpdokposdkfpsdokfpsokshfigushfzdiziyuezaiozo} entails condition $(ii)$. Finally, one has the equality:
\begin{equation}
\label{fpkerfpekrfpo}
 I:=\int^{{(2e)}^{-1}}_{0}  R^{{\scriptscriptstyle-1/2}}_{u}\ dR_{u} = \int^{{(2e)}^{-1}}_{0} \ R^{{\scriptscriptstyle-1/2}}_{u}\cdot R'_{u} \ du = 2\ \int^{{(2e)}^{-1}}_{0} \ (R^{{\scriptscriptstyle1/2}}_{u})' \ du = 2\cdot R^{{\scriptscriptstyle1/2}}_{1/2e}.
\end{equation}

\section{Tanaka formula \& Local times}
\label{Tanaka}
The main results of this section are a Tanaka formula (Theorem 
\ref{formule de tanaka}) and two occupation time formulas (Theorem \ref{oedzemlijz}).
In Subsection \ref{tanakaformula} 
we perform a complete comparison between our Tanaka formula and the ones provided so far, for general Gaussian processes (namely \cite[Theorem $4.3$]{LN12} and, to a lesser extent, in \cite[Theorem $36$]{MV05}). In particular, it will be shown that the Tanaka formula we provide here not only fully extends the one provided in \cite[Theorem $36$]
{MV05} but is also less restrictive, in many places, than the one provided by \cite[Theorem $4.3$]{LN12}. We then apply our Tanaka formula to 
the Gaussian processes of ``reference'' and compare the resulting 
Tanaka formulas hence obtained with the Tanaka formulas established so far, for 
these particular Gaussian processes.
In Subsection  \ref{ozdeiojezoijziojéazpokdzlkn} we 
provide a definition and a basic study of both (non weighted) and 
weighted local times of elements $G$ in $\mathscr{G}$. 
After having provided an integral representation of each of 
these two processes,  
we establish, under classical assumptions, the belonging of both these 
local times, as two parameters processes, to $L^{2}(\R\times
\Omega, \lambda \otimes \mu)$. The proof, that (non weighted) local 
time of $G$ belongs to 
$L^{2}(\R\times\Omega, \lambda \otimes \mu)$ we give here, is not only 
new and based on White noise analysis; it also offers the advantage to 
be easily adaptable to the weighted local time of $G$. It is then 
shown in what extent the results we provide in this Subsection are new, 
and we also give a complete comparison of our results to the one 
obtained in \cite{LN12}, for weighted local times (that are the 
only local times considered in this latter reference). Finally, we apply our result to the Gaussian processes of ``reference'' and show how we 
extend the results known for them, so far.
\vspace{0.5ex}

\hspace{0.5cm}In addition to the references \cite{GH,MARO} and  
\cite[Chapter VI]{RY}, that provide a complete overview on local times, a summary on fractional Brownian motion's 
local time (that includes the concept of weighted 
fractional local time) can be found in \cite[Section $7.8$]
{COU07} (see also references therein, in particular \cite{CNT01}). Note that a 
study of fBm's local time in the 
white noise theory's framework can be found in 
\cite{HOS05}.
The non-weighted local times of particular mBms have been studied in 
\cite{BDG06,BDG07} in the case of an harmonizable mBm while the 
weighted local time of a Volterra-type mBm has been studied in \cite{dozzi}. 
This approach has been extended to multifractional Brownian sheets in 
\cite{MWX08}.


%
%
%
%
%
%
%
%
%
\smallskip

Through this section, let $T$ denotes a positive real and recall that $\z:= \Z_{R} \cap [0,T]$.

\subsection{Tanaka Formula}
\label{tanakaformula}

A consequence of Theorem \ref{zpzpkpskpsokpdokposdkfpsdokfpsokshfigushfzdiziyuezaiozo} is the following Tanaka formula.

\begin{theo}[Tanaka formula for $G$]
\label{formule de tanaka}
Let $T>0$ be such that $[0,T]\subset \rR$ and $c$ be real number. Assume that Assumption \Aa holds and that the map $t\mapsto R_{t}$ is both continuous and of bounded variations on $[0,T]$ and such that: 
\vspace{-1ex}

\bit
\iti $\lambda(\z) = \alpha_{R}(\z)=0$,
\vspace{-1ex}
\itii $t\mapsto {R^{\scriptscriptstyle-1/2}_{t}} \in L^1([0,T], dR_t)$,
\vspace{-0.5ex}
\itiii $\exists \ q \in \N$ such that $t\mapsto {|g'_t|}_{-q} \ R^{-1/2}_{t} \in L^1([0,T],dt)$.
\eit
 \vspace{-1ex}
 
Then, the following equality holds in $(L^2)$:
\vspace{-1ex}
\begin{equation}
\label{tanfor}
|G_{t} - c| = |c| + \int^{T}_{0} \   \text{sign}(G_{t} - c) \  \di G_{t} +  \int^{T}_{0}  \  \delta_{\{c\}}(G_{t}) \   dR_{t},
\end{equation}
where the function sign is defined on $\R$ by $\si(x) :=\i1_{\R^*_{+}}(x) - \i1_{\R_{-}}(x)$.
\end{theo}

\begin{pr}
This is a direct application of Theorem \ref{zpzpkpskpsokpdokposdkfpsdokfpsokshfigushfzdiziyuezaiozo} with $F:x \mapsto |x-c|$.
\end{pr}

\begin{rem}
\label{fierjoerij}
Thanks to \eqref{fpkerfpekrfpo}, we clearly see that 
condition $(ii)$ is satisfied if $t\mapsto R_t$ is absolutely continuous.
\end{rem}

\vspace{-3ex}

%
%
%

\subsubsection*{Comparison with Tanaka formula for general and particular Gaussian processes}

As Equality \eqref{zpedokzpdokzozek} will show, one can write \eqref{tanfor} with the (weighted) local time instead of the last term on the right hand side. 

Because the approach used in \cite{LN12} is more intrinsic (since all the assumptions are made exclusively on the variance 
function $R$), it is not easy to compare conditions (H1), (H3a) \& (H4), 
required by \cite[Theorem 4.3]{LN12}, to conditions (i) to (iii) of 
Theorem \ref{formule de tanaka}. However, one can note the following facts. It is clear that one does not 
need here to make any assumption on the covariance function $(t,s)
\mapsto R_{t,s}$ but only on the variance function $t\mapsto R_{t}$. In \cite{LN12} $(t,s)\mapsto R_{t,s}$ is assumed to be continuous on ${[0,T]}^{2}$, while only $t\mapsto R_{t}$ is assumed to be continuous on ${[0,T]}$ here. 
Moreover, we do not need to assume that $t\mapsto R_{t}$ is increasing, 
nor such that $\z=\{0\}$. Besides, 
the fact that $t\mapsto R_{t}$ is of bounded variations on $[0,T]$ is 
sufficient here and one thus does not need to make Assumption (H3a). Finally we do not need any assumption of the type of \cite[(H4)]{LN12} to get a Tanaka formula.

%
%
%
%
%
%
%
%
%
%
Besides, in view of results of \cite[Theorem 5.6]{JL17a} 
and \eqref{fpkerfpekrfpo}, it is clear 
that Theorem \ref{tanakaformula} applies to all $\mathscr{V}_{\gamma}$ - processes and, therefore, that Equality \eqref{tanfor} (or if one prefers Equality \eqref{zpedokzpdokzozek} below) fully generalizes 
 \cite[Theorem $36$]{MV05} (with $F:x \mapsto {(x-c)}_{+}$ instead of $x \mapsto {|x-c|}$).
 
Note moreover that Theorem \ref{tanakaformula} also generalizes the classical Tanaka formula 
for Brownian motion (see \cite[Theorem $1.5$]{RY}), the one given for 
fBm in \cite[Corollary $4.8$]{ben1} for any Hurst index $H$ in $(0,1)$, the one provided in \cite[Theorem 
$12$]{dozzi} for Volterra type mBm (it is assumed in there that $h$ takes 
values in $(1/2,1)$) and  the one given in \cite[Theorem $6.1$]{JLJLV1} 
for normalized mBm.

\subsection{Local times of Gaussian processes in $
\mathscr{G}$}
\label{ozdeiojezoijziojéazpokdzlkn}

In all this subsection one does not assume anymore that Assumptions \Aa holds. $G$ is just assumed to be of the form \eqref{erofkorekp} and such that the map $g:t\mapsto g_{t}$, defined at the beginning of Subsection \ref{sjjjj},
 is continuous on $\rR$.  One assumes, through this section, that  $[0,T]\subset \rR$. Let us first define the local time of $G$.
\vspace{-1ex}



%
%

%
%

%

%

\begin{defi}{(non weighted) local time of $G$}
\label{lt}

The (non weighted) local time of $G$ at any point $a\in\R$, up to time $T$, denoted $\ell^{(G)}_{T}(a)$, is defined by:
\begin{equation*}
\ell^{(G)}_{T}(a):=\underset{\varepsilon\to 0^{+}}{\lim} \frac{1}{2\varepsilon} \mathcal{\lambda}(\{ s\in [0,T]; \  G_{s} \in (a-\varepsilon, a+\varepsilon)  \}), 
\end{equation*}
where $\lambda$ denotes the Lebesgue measure on $\R$ and where the limit holds in \oS\hspace{-0.1cm}, when it exists.
\end{defi}


%
%


Having in mind the definition of generalized functional of $G$ (given in 
Theorem-Definition \ref{lmqsdlmflkjdg}) one has the following result.



%
%
%

\begin{prop}{[Integral representation on $[0,T]$ of non weighted local time]}
\label{ozrisjoerzijo} 

Assume that the variance map $s\mapsto R_{s}$ is continuous on $[0,T]$ and such that the Lebesgue measure of $\z$ is equal to $0$, then:
\vspace{-0.25cm}

\ben
\item The map  $s\mapsto \delta_{a}(G_{s})$ is \oS-integrable on $[0,T]$ for every $a\in \R^{*}$. The map $s\mapsto \delta_{0}(G_{s})$ is \oS-integrable on $[0,T]$ if $s\mapsto {R^{\scriptscriptstyle -1/2}_{s}}$ belongs to $L^{1}([0,T],ds)$.  
\item The following equality holds in \oS, 
for every $a\in \R^{*}$:
\begin{equation}
\label{oeirjeorifj}
\ell^{(G)}_{T}(a) = \int^{T}_{0} \delta_{a}(G_{s}) \ ds,
\end{equation}
\een
\vspace{-2ex}
Equality \eqref{oeirjeorifj} still holds 
 for $a=0$ if $s\mapsto {R^{\scriptscriptstyle-1/2}_{s}} \in L^{1}([0,T],ds)$.
\end{prop}

\begin{rem}
The reason why we have to assume that the Lebesgue measure of $\z$ is equal 
to $0$ lies in the fact that, when $F$ belongs to $\sS'(\R)$, $F(G_{t})$ is not defined if $t$ belongs to $\z$ (Theorem-Definition \ref{lmqsdlmflkjdg}).
However, for any $a$ in $\R^{*}$\hspace{-0.1cm},  if $\delta_{a}(G_{s})$ has 
a limit in \oS, as $s$ tends to $s_{0}\in \z$, one may extend $\delta_{a}
(G_{s})$ at point $s_{0}$ by setting $\delta_{a}(G_{s_{0}}) := 0$. We can not say anything in general, when $a=0$.

\end{rem}
%
Before giving the proof of Proposition \ref{ozrisjoerzijo}, let us state the following result, the proof of which is obvious and therefore left to the reader, since $t\mapsto R_{t}$ is a continuous function.

\begin{lem}
\label{eofijeofjeroifjerofjer}
One has the following result:
\vspace{-4.5ex}

\begin{align}
\label{epoferooerijeofoe}
 \forall\ \eta \in \sS(\R),\ \forall a \in \R^{*},\ &\forall \varepsilon \in (0,\tfrac{|a|}
{2}),\ \exists \ \gamma \in (0,\varepsilon), \  \text{such that: }\ \forall (s,y) \in D(\z,\gamma)\times D(a,\gamma),  \notag \\
& -{(y-<g_{s},\eta>)}^{2} < -a^{2}/16,
 \end{align}

where we have set $D(a,\gamma):=(a-\gamma, a+\gamma)$, $D(\z,\gamma):=\{u\in [0,T]; \ d(u,\Z^{T}_R) < \gamma\}$ and where $d(A,x)$ denote the distance between the subset $A$ of $\R$ and the real $x$.
\end{lem}

{\bfseries Proof of Proposition \ref{ozrisjoerzijo}.}

{\bfseries 1.} Thanks to \cite[Theorem $7.3$]{Kuo2}, we know that $s\mapsto 
\delta_{a}(G_{s})$ is an \oS-process, the $S$-transform of which is equal, for every $(a,\eta,s) \in \R\times\sS(\R)\times \w$, to:
\begin{equation}
\label{oeivjeorijv}
S(\delta_{a}(G_{s}))(\eta)= 
\tfrac{1}{\sqrt{2\pi R_{s}}} \exp\big\{\tfrac{-{(a-<g_{s},\eta>)}^{2}}{2 R_{s}}\big\},
\end{equation}
where $\w:=\{u\in [0,T];\ R_{u}>0\}$. When $a=0$, the right hand side of \eqref{oeivjeorijv} is bounded by $R^{\scriptscriptstyle-1/2}_{s}$, which belongs, by assumption, to $L^{1}([0,T], dt)$.
 Theorem \ref{peodcpdsokcpodfckposkcdpqkoq} then applies and establishes {\bfseries 1}. in the case where $a=0$.

Let $a$ be in $\R^{*}$. The weak measurability of $t\mapsto \delta_{a}(G_{t})$ is clear. In view of \cite[p. $246$]{Kuo2}, one just has to show that $t\mapsto S(\delta_{a}(G_{t}))(\eta)$ belongs to $L^{1}([0,T],dt)$, for every $\eta$ in $\sS(\R)$. Let $\eta$ be in $\sS(\R)$ and $\varepsilon$  be in $(0,\tfrac{|a|}{2})$. In view of \eqref{oeivjeorijv} and thanks to Lemma \ref{eofijeofjeroifjerofjer}, one knows that, for every $s_{0}$ in $\z$, there exists a neighborhood $\mathcal{V}_{s_{0}}$  of $s_{0}$ on which $|S(\delta_{a}(G_{s}))(\eta)| \leq  R^{\scriptscriptstyle-1/2}_{s} \ e^{-a^{2} / 32R_{s}}=:l_{a,\eta}(s)$.
Since $l_{a,\eta}$ belongs to $L^{1}([0,T], dt)$ for every $\eta$ in $\sS(\R)$, one has established  {\bfseries 1}. in the case $a\neq 0$ and therefore ends the proof of {\bfseries 1}.
\smallskip

{\bfseries 2.} Let $\varepsilon >0$, and $a$ be a real number. Define 
\vspace{-0.5cm}

\begin{align*}
I_{\varepsilon}(a) \ :=\lambda(\{ s\in [0,T]; \ G_{s} \in D(a,\varepsilon)\});\ \hspace{0.2cm}\text{and }\ \hspace{0.2cm}  \Phi_{\varepsilon}(a) \stackrel{\text{\textcolor{white}{a}\oS}}{: =}  \int_{\R}\ \big( \int^{T}_{0}\   \i1_{(a-\varepsilon, a+\varepsilon)} (y)\ \delta_{y}(G_{s}) \ ds \big) \ dy. 
\end{align*}

%

A simple application of Theorem \ref{peodcpdsokcpodfckposkcdpqkoq} and 
Fubini's theorem shows that $\Phi_{\varepsilon}(a)$ is well defined, for every $a\in\R$. Besides, since $I_{\varepsilon}(a)= \int^{T}_{0}\  \i1_{D(a,\varepsilon)}(G_{s}) \ ds$, one has:

\vspace{-0.5cm}
\begin{align}
\label{oeoifjerojeroi}
I_{\varepsilon}(a)&= \int^{T}_{0}\  \left(\int_{\R}\ \i1_{(a-\varepsilon, a+\varepsilon)} (y)\ \delta_{y}(G_{s}) \ dy \right)\ ds,
\end{align}

where we used the following lemma, the proof of which can be found in Appendix \ref{appendiceS}.

\begin{lem}
\label{eporkerpok}
Let $j \in\sS'(\R)$ be of function type and let $g$ be in $L^{2}(\R)$. Define $\Psi_{j}:\R\rightarrow \text{\oS}$ by setting $\Psi_{j}(y)= j(y)\ \delta_{y}(<.,g>)$. If $\Psi_{j}$ is \oS-integrable on $\R$, then one has the equality
\vspace{-0.15cm}

\begin{equation}
\label{peoreorfeorij}
 \int_{\R} \ j(y)\ \delta_{y}(<.,g>)\ dy \stackrel{\text{\oS}}{=} j(<.,g>).
\end{equation}

Moreover, if  $j(<.,g>)$ belongs to $(L^{2})$, then \eqref{peoreorfeorij} is also true in $(L^{2})$.
\end{lem}

 A direct application of Fubini's Theorem and \eqref{oeoifjerojeroi} shows the equality of the $S$-transform of $\Phi_{\varepsilon}(a)$ and $I_{\varepsilon}(a)$. The equality in \oS results of the injectivety of the $S$-transform. It then remains to show that $\ds 
\underset{\varepsilon\to 0}{\lim} \tfrac{I_{\varepsilon}(a)}{2 \varepsilon}$ exists in \oS. Besides, let $(\eta, a)$ be in $\sS(\R)\times\R$. Let $\varepsilon$ be in $(0,1)$ such that $D(a,\varepsilon)\subset \R^{*}_{+}$ or $D(a,\varepsilon)\subset \R^{*}_{-}$, if $a\in \R^{*}$. Using $(ii)$ of Lemma 
\ref{dkdskcsdckksdksdmksdmlkskdm} and \eqref{oeivjeorijv} we get:

\vspace{-0.5cm}

\begin{align}
\label{lamballe}
S(\tfrac{I_{\varepsilon}(a)}{2 \varepsilon})(\eta)&= \frac{1}{2\varepsilon}  \hspace{-0.5ex}\int_{D(a,
\varepsilon)}  \hspace{-1ex}\bigg(\int^{T}_{0}  \hspace{-1ex} \tfrac{1}{\sqrt{2\pi R_{s}}} \exp\{\tfrac{-1}{2 
R_{s}} {(y-<g_{s},\eta>)}^{2}\}   \ ds\bigg)\ dy =: \frac{1}{2\varepsilon}  \hspace{-0.5ex} \int_{D(a,
\varepsilon)} \hspace{-1ex} \big(\int^{T}_{0}\hspace{-1ex} \rho_{\eta}(y,s)\ ds \big) dy.
\notag\\
\end{align}

 \vspace{-1ex}
 
We will use the following lemma, the proof of which can be found in Appendix \ref{appendiceS}.

\begin{lem}
\label{pzokpo}
For every $\eta$ in $\sS(\R)$, define the map $G_{\eta}:\R^{*} \rightarrow \R$, by setting 
$G_{\eta}(y):=\int^{T}_{0} \rho_{\eta}(y,s)\ ds$. Then $G_{\eta}$ is continuous on $\R^{*}$. Moreover, if $s\mapsto 
{R^{\scriptscriptstyle-1/2}_{s}}$ belongs to $L^{1}([0,T]\backslash\z,ds)$, then one can extend $G_{\eta}$ to $\R$ by setting $G_{\eta}(0):= \int^{T}_{0} \rho_{\eta}(0,s)\ ds$. Furthermore, this extension is continuous on $\R$.
\end{lem}

The equality $G_{\eta}(a)=S(\int^{T}_{0} \delta_{a}(G_{s}) \ ds)(\eta)$ and the continuity of $G_{\eta}$, on $D(a,\varepsilon)$, allows us to write $ 
\underset{\varepsilon\to 0}{\lim}\  \frac{1}{2\varepsilon} \int_{D(a,
\varepsilon)}  G_{\eta}(y)\ dy=S(\int^{T}_{0} \delta_{a}(G_{s}) \ ds)(\eta)$. Thanks to \eqref{lamballe}, one finally gets $ 
\underset{\varepsilon\to 0}{\lim}\  S(\tfrac{I_{\varepsilon}(a)}{2 \varepsilon})(\eta)=S(\int^{T}_{0} \delta_{a}(G_{s}) \ ds)(\eta)$. Moreover, using the same computations as in the proof of Lemma \ref{pzokpo}, it is easy to show that, for every $a$ in $\R^{*}$, 
\vspace{-2ex}

\begin{equation*}
\forall \eta \in \sS(\R), \ \forall y \in D(a,|a|/2), \hspace{0.5cm}  |G_{\eta}(y)| \leq M_{a}\ e^{{|\eta|}_{0}},
\end{equation*}

where $M_{a}:= T \ (8(1+e){|a|}^{-1} + \sup\{{R^{\scriptscriptstyle-1/2}_{s}} {\exp}(-a^{2}/32 R_{s}) \}) $. One then gets: $|S(\tfrac{I_{\varepsilon}(a)}{2 \varepsilon})(\eta)|\leq M_{a}\ e^{{|\eta|}_{0}}$.
%
%
\cite[\hspace{-0.05ex}Theorem $8.6$]{Kuo2} then allows us to write {\small$
\underset{\varepsilon\to 0}{\lim}\  \tfrac{I_{\varepsilon}(a)}{2 \varepsilon} 
\stackrel{\text{\oS}}{=} \int^{T}_{0} \delta_{a}(G_{s})\hspace{0.1ex} ds$}, which ends 
the proof.\hfill { $\square$}

\begin{rem}
We deduce from the proof of Proposition \ref{ozrisjoerzijo} that $\ell^{(G)}
_{T}(a)$ exists for every $a$ in $\R^{*}$; and also for $a=0$, if $s\mapsto 
{R^{\scriptscriptstyle-1/2}_{s}}$ belongs to $L^{1}([0,T],ds)$. Besides, the 
example of the Gaussian process $(G_{s})_{s\in\R}$, defined by setting 
$g_{s}:=\sum_{k\geq 0} {s} {(1+k)}^{-1}\ e_{k}$, shows that the condition $s
\mapsto {R^{\scriptscriptstyle-1/2}_{s}}$ belongs to $L^{1}([0,T],ds)$ is required in order 
that the map $s\mapsto \delta_{0}(G_{s})$ is \oS-integrable on $[0,T]$.
\end{rem}

\begin{theodef}{[Weighted local time of $G$ and its integral representation in \oS]}
\label{ltw} 

Assume that the map $s\mapsto R_{s}$ is continuous and of bounded variations on $[0,T]$ and such that $\alpha_{R}(\z)$ is equal to $0$. Then:
\vspace{-1ex}
\ben
\item The map  $s\mapsto \delta_{a}(G_{s})$ is \oS-
integrable on $[0,T]$ with respect to the measure $\alpha_{R}$, for every $a\in\R^{*}$. The map $s\mapsto \delta_{0}
(G_{s})$ is \oS-integrable on $[0,T]$  with respect to the measure  $\alpha_{R}$, if $s\mapsto  R^{\scriptscriptstyle-1/2}_{s}$ belongs to 
$L^{1}([0,T], dR_{s})$.

\item For every $a\in\R$, when  $s\mapsto  \delta_{a}(G_{s})$ is \oS-integrable on $[0,T]$  with respect to the measure  $\alpha_{R}$, one can define the weighted local time of $G$ at point $a$, up to time $T$, denoted  $\mathscr{L}^{(G)}_{T}(a)$, as being the \oS- process defined by setting:
\een
\vspace{-3ex}

\begin{equation}
\label{oeirdeezjeorifj}
\mathscr{L}^{(G)}_{T}(a) := \int^{T}_{0}  \ \delta_{a}(G_{s}) \ dR_{s},
\end{equation}
\vspace{-2ex}

where the equality holds in \oS.

\end{theodef}

\begin{pr}
Since the proof of Point $1$ follows exactly the same steps as proof of Proposition \ref{ozrisjoerzijo} (one just has to replace the Lebesgue measure by $\alpha_{R}$), it is left to the reader.
\end{pr}
Assume that $s\mapsto R_{s}$ is absolutely continuous on $[0,T]$ and such that its derivative, denoted $R'$, is differentiable \textit{a.e.} on $[0,T]$.  If the following conditions are fulfilled:
\vspace{-2ex}
\bit
\iti $s\mapsto R'_{s}\ \delta_{a}(G_{s})$ is \oS-integrable on $[0,T]$,
\itii $s\mapsto R''_{s}\ \ell^{(G)}_{s}(a)$ is \oS-integrable on $[0,T]$,
\itiii $s\mapsto R'_{s}\ \ell^{(G)}_{s}(a)$ has a limit in $0$ and $T$ (denoted respectively {$\scriptstyle \ell^{(G)}_{s}(a)  \ R'_{s}\big|_{s=0}$ and  $\scriptstyle \ell^{(G)}_{s}(a)  \ R'_{s}\big|_{s=T}$}),
\eit

then an integration by parts yields:
\vspace{-0.15cm}

\begin{equation}
\label{oeirjodddd}
 \mathscr{L}^{(G)}_{T}(a) = \ell^{(G)}_{s}(a)  \  R'_{s}\big|_{s=T}  -  \ell^{(G)}_{s}(a) R'_{s} \big|_{s=0}
- \int^{T}_{0} \ \ell^{(G)}_{s}(a)\ {R''_{s}}  \ ds,
\end{equation}

where $R''$ denotes the second derivative of $R$.

\begin{rem}
\label{pokopkpokopkdzyytzedfzueudpéodheug}
{\bfseries 1.} In view of the definition of weighted local time, it is clear that \eqref{tanfor} now reads:
\vspace{-2ex}

\begin{equation}
\label{zpedokzpdokzozek}
|G_{t} - c| = |c| + \int^{T}_{0} \   \text{sign}(G_{t} - c) \  \di G_{t} +   \mathscr{L}^{(G)}_{T}(c).
\end{equation}

As we stated at the end of Subsection \ref{tanakaformula}, one recovers hence the Tanaka formula provided in \cite[Theorem $36$]{MV05} and in \cite[Theorem 4.3]{LN12}.

{\bfseries 2.} The definition of weighted local time of $G$ (Equality \eqref{oeirdeezjeorifj}) generalizes the one given, for fBm, in \cite[Definition $14$]{COU07}. Note moreover that Equality \eqref{oeirjodddd} above generalizes, to all Gaussian processes which belong to $
\mathscr{G}$, the relation between weighted and non-weighted 
local times given, for fBm, in \cite[Definition $14$]{COU07}.
\end{rem}
%
%


\section{Occupation times formulas and Comparison to previous results}
\label{occupationetcomparaison}
\subsection{Occupation times formulas}
\label{occupation}
The following theorem establishes that, under suitable assumptions, both weighted and non-weighted local times are $(L^{2})$ random variables. It also provides an occupation time formula for both $\ell^{(G)}$ and $\mathscr{L}^{(G)}$. Denote $\mathscr{M}_{b}(\R)$ the set of positive Borel functions defined on $\R$ and recall that $\lambda$ denotes the Lebesgue measure on $\R$.

\begin{theo}
\label{oedzemlijz}
Assume that the function $s\mapsto R_{s}$ is continuous on [0,T].
\vspace{-2ex}
\bit      
\iti Assume that $\lambda_{R}(\z)=0$ and that $s\mapsto {R^{\scriptscriptstyle-1/2}_{s}}\in L^{1}([0,T],ds)$.
 If $\E\big[\int_{\R} {\big| \int^{T}_{0}  e^{i
\xi G_{s}}   ds\big|}^{2} d\xi \big]\hspace{-0.75ex}<+\infty$, then the map $a\mapsto \ell^{(G)}_{T}(a)$  belongs to $L^{2}
(\lambda \otimes \mu)$. 
Moreover one has the following equality, valid for $\mu$-a.e. $\omega$ in $
\Omega$,
\begin{equation}
\label{zpckzpkzepokzepoczepock}
 \forall\ \Phi\ \in \mathscr{M}_{b}(\R), \ \int^{T}_{0}\  \Phi(G_{s}(\omega)) \ ds =  \int_{\R}\  \ell^{(G)}_{T}(y)(\omega) \cdot \Phi(y) \ dy.
\end{equation}

\vspace{-0.2cm}
\itii Assume that $\alpha_{R}(\z)=0$. Assume moreover that $s\mapsto R_{s}$ is of bounded variations on $[0,T]$, and 
such that $s\mapsto {R}^{\scriptscriptstyle-1/2}_{s} \in L^{1}([0,T], dR_{s})$.  
 If $
\E\big[\int_{\R}\ {\big| \int^{T}_{0}  \ e^{i\xi G_{s}}  \ dR_{s} \big|}^{2}\ d\xi 
\big] <+\infty$,  then $a\mapsto\mathscr{L}^{(G)}_{T}(a)$  belongs to $L^{2}
(\lambda \otimes \mu)$. Moreover one has the following equality, valid for $
\mu$-a.e. $\omega$ in $
\Omega$,
\begin{equation}
\label{eovjeojoijovjeroivejro}   \forall \ \Phi\in \mathscr{M}_{b}(\R), \ \int^{T}_{0}\  \Phi(G_{s}(\omega)) \ dR_{s} =  \int_{\R}\  \mathscr{L}^{(G)}_{T}(y)(\omega) \cdot \Phi(y) \ dy.
\end{equation}
\eit
\end{theo}

%

\begin{rem}
\label{pokerpofkerpoferpoferpofkerp}
${\text{{\bfseries{1.}}}}$\ The first consequence of Theorem \ref{oedzemlijz} is that both $\ell^{(G)}_{T}(a)$ and $\mathscr{L}^{(G)}_{T}(a)$ are $(L^{2})$ random variables for every real $a$. In a forthcoming paper we discuss the joint continuity of processes
$(T,a)\mapsto \ell^{(G)}_{T}(a)$ and $(T,a)\mapsto \mathscr{L}^{(G)}_{T}(a)
$. In particular, under assumptions of local non-determinism\footnote{see \cite[\S3.24]{GH} and references therein for an overview of this 
notion}, it is shown that both map $a\mapsto \ell^{(G)}_{T}
(a)$ and $a\mapsto\mathscr{L}^{(G)}_{T}(a)$ are continuous (the 
continuity at point $0$ requires that $t\mapsto {R}^{-1/2}_{t}$ belongs to 
$L^{1}([0,T], dt)$ or to  $L^{1}([0,T], dR_{t})$).

${\text{{\bfseries{2.}}}}$ Since $\E[X]=S(X)(0)$, for every $(L^{2})$ random variable, it is easy to get, for almost every real $a$, the equalities:
\vspace{-3.5ex}

\begin{align*}
\E[\ell^{(G)}_{T}(a)] &= \int^{T}_{0}  \tfrac{1}{\sqrt{2\pi\ R_{s}}} \ \exp\{-\tfrac{a^{2}}{2\  R_{s}}   \} \ ds; \hspace{3.5ex} \& \hspace{3.5ex} \E[\mathscr{L}^{(G)}_{T}(a)] &=   \int^{T}_{0}  \tfrac{1}{\sqrt{2\pi\ R_{s}}}\ \exp\{-\tfrac{a^{2}}{2\  R_{s}}   \} \ dR_{s},
\end{align*}
when assumptions of Theorem \ref{oedzemlijz} are fulfilled,

${\text{{\bfseries{3.}}}}$\ 
As we stated right after the proof of Theorem-Definition \ref{ltw},  the 
map  $t\mapsto {R^{\scriptscriptstyle-1/2}_{t}}$ belongs to $L^{1}([0,T], dR_{t})$
 if $t\mapsto R_{t}$ is absolutely continuous. In 
this latter case \eqref{eovjeojoijovjeroivejro} holds under the single 
condition $
\E\big[\int_{\R}\hspace{-0.5ex} {\big| \hspace{-0.5ex}\int^{T}_{0}   e^{i\xi 
G_{s}}   dR_{s} \big|}^{2} d\xi 
\big]\hspace{-0.5ex}<\hspace{-0.5ex}+\infty$, which is equivalent to $\int^{T}_{0} \int^{T}
_{0} {(\Delta(t,s))}^{\scriptscriptstyle-1/2}  dR_{s} \ dR_{t}\hspace{-0.75ex}<+\infty$, where $\Delta(t,s):= 
\E[{(G_{t}-G_{s})}^{2}]$. As it is shown in 
\cite[(21.9) \& (21.10)]{GH} any of these two 
latter conditions is necessary and sufficient to ensure the existence of a 
weighted local time, that is an $(L^{2})$ random variable. This shows that 
the assumptions we made on the variance function $R$ are almost 
minimal.

${\text{{\bfseries{4.}}}}$\ 
The condition $t\mapsto {R^{\scriptscriptstyle-1/2}_{t}}\in L^{1}([0,T], dt) \cap  L^{1}([0,T], dR_{t})$
  is verified for all the Gaussian processes of ``reference'' (see Proposition \ref{pofkeporfkerpofk} below for the case of $\mathscr{V}_{\gamma}$ - processes).
%
%
\end{rem}



Let us give the proof of Theorem \ref{oedzemlijz} before discussing the results provided therein and before making the comparison with the results obtained in \cite[\S 4]{LN12}.


\begin{pr}
Let us start with the proof of $(i)$. Define, for every $(\omega,a)$ in $\Omega\times \R$,  the map denoted $f^{(\omega)}_{T}: 
\R\rightarrow \R$ and, for every positive integer $n$, the family of maps denoted $f^{(\omega)}_{T,n}: \R\rightarrow 
\R$ by setting: $f^{(\omega)}_{T}(\xi):= \frac{1}{2\pi}
\int^{T}_{0} e^{i  \xi G_{s}(\omega)} \ ds$ and 
$f^{(\omega)}_{T,n}(\xi):= e^{-{\xi}^{2}/n}\ f^{(\omega)}_{T}
(\xi)$. Recalling the notations about Fourier transform defined in Section \ref{cddcwcdsd}, define also:
\vspace{-0.5cm}

\begin{align}
K^{(\omega)}_{T}(a):= \cF\big(f^{(\omega)}_{T}\big)(a) ; \hspace{0.5cm}\&  \hspace{0.5cm} K^{(\omega)}_{T,n}(a):= \cF\big(f^{(\omega)}_{T,n}\big)(a).
\end{align}

It is clear that $\underset{n\to+\infty}{\lim} {K^{(\omega)}
_{T,n}} = K^{(\omega)}_{T}$, in ${L^{2}(\R)}$, since $\underset{n\to+\infty}{\lim^{(L^{2)}}} {f^{(\omega)}_{T,n}} = f^{(\omega)}_{T}$.
Denote $K_{T,n}(a)$ 
(\textit{resp.} $K_{T}(a)$) the random variable $\omega \mapsto 
K^{(\omega)}_{T,n}(a)$ (\textit{resp.}  $\omega \mapsto 
K^{(\omega)}_{T}(a)$).

{\bfseries Step 1:} $\underset{n\to+\infty}{\lim} {K_{T,n}}(a) =  \ell^{(G)}_{T}(a)$, where the equality holds in \oS.

Let $\omega$ in $\Omega$. Since $f^{(\omega)}_{T,n}$ belongs to $L^{1}(\R)\cap L^{2}(\R)$, one can write ${K^{(\omega)}_{T,n}}(a) = \int_{\R} f^{(\omega)}_{T,n}(\xi) \ e^{-ia\xi}\ d\xi$. For every $\eta$ in $\sS(\R)$, using Lemma \ref{dkdskcsdckksdksdmksdmlkskdm}, we then have:
\vspace{-0.3cm}

\begin{align*}
S(K_{T,n}(a))(\eta) = \int_{\R}\  \frac{1}{2\pi} e^{-ia\xi -{\xi}^{2}/n}  \bigg(\int^{T}_{0}\  S(e^{i\xi G_{s}})(\eta) \ ds \bigg)\ d\xi.
\end{align*}

Using Equality \eqref{oeij12} and Fubini's theorem, that both obviously apply here, we get:

\vspace{-2ex}

\begin{align*}
S(K_{T,n}(a))(\eta) = \tfrac{1}{\sqrt{2\pi}} \int^{T}_{0}  \bigg(\int_{\R} \tfrac{1}{\sqrt{2\pi}} e^{-\frac{1}{2}\big(\xi^{2}(\frac{2}{n}+R_{s}) + 2i \xi (a-<\eta,g_{s}>)  \big)}  d\xi\bigg)  ds =: \tfrac{1}{\sqrt{2\pi}} \int^{T}_{0}  J_{n,\eta}(s) \ ds.
\end{align*}

Using the change of variable $u:= \xi\cdot \sqrt{\varepsilon_{n}+R_{s}}$, where $\varepsilon_{n}:=\frac{2}{n}$, in $J_{n,\eta}$ yields:
\vspace{-3ex}

\begin{align*}
S(K_{T,n}(a))(\eta) = \tfrac{1}{\sqrt{2\pi}} \int^{T}_{0}\   \tfrac{{\exp}\{-\tfrac{ {(a-<\eta,g_{s}>)}^{2}  }{2 (\varepsilon_{n}+R_{s}) }\}}{\sqrt{\varepsilon_{n}+R_{s}} }  \ ds.
\end{align*}


Using the same arguments that the one we used in the proof 
of Lemma \ref{pzokpo} 
allows us to apply the Lebesgue dominated convergence. 
We therefore get, using \eqref{oeivjeorijv} and then \eqref{oeirjeorifj},

\vspace{-0.35cm}

\begin{align}
\label{zoijfeorijeroiceroij}
\underset{n\to+\infty}{\lim} S(K_{T,n}(a))(\eta) &=  
\int^{T}_{0}\ \frac{1}{\sqrt{2\pi R_{s}}}     \exp \{  -\tfrac{ {(a-<
\eta,g_{s}>)}^{2}  }{2 R_{s} } \}
 \ ds
 =S(\int^{T}_{0} \delta_{a}(G_{s}) \ ds)(\eta),\notag\\
 &=S(\ell^{(G)}_{T}(a))(\eta).
\end{align}

The domination used in the proof of Lemma \ref{pzokpo}
 applies here also and prove the existence of a function $L$, which belongs to in $L^{1}([0,T])$, such that:

\vspace{-0.35cm}

\begin{align}
\label{zoijfeorijeroiceroij2} 
\forall(p,n,\eta,a)\in \N\times \N^{*}\times \sS(\R)\times \R,\hspace{0.5cm}  |S(K_{T,n}(a))(\eta)| \leq \ e^{{|\eta|}^{2}_{p}}\cdot \int^{T}_{0} L(u)\ du.
\end{align}

Formulas \eqref{zoijfeorijeroiceroij} and \eqref{zoijfeorijeroiceroij2} allow us finally to apply \cite[Theorem $8.6.$]{Kuo2} and we then establish that  $\underset{n\to+\infty}{\lim} {K_{T,n}}(a) =  \ell^{(G)}_{T}(a)$, in \oS.
\vspace{0.1cm}

{\bfseries Step 2:} Let us show that there exists a strictly increasing function $\varphi:\N^{*}\rightarrow\N^{*}$ such that, for almost every real $a$, $\underset{n\to+\infty}{\lim} {K_{T,\varphi(n)}}(a) = {K_{T}}(a)$, in $(L^{2})$.  

From the first step we know that, for $\mu$ - almost every $\omega$, $\underset{n\to+\infty}{\lim} {\|{K^{(\omega)}
_{T,n}} - K^{(\omega)}_{T}\|}^{2}_{L^{2}(\R)} =0$. From the other hand, one gets $\underset{n\in\N^{*}}{\sup} {\|{K^{(\omega)}
_{T,n}} \|}_{L^{2}(\R)} \leq {\|{K^{(\omega)}
_{T}} \|}^{2}_{L^{2}(\R)}$.
Since, by assumption,  $\omega \mapsto {\|{K^{(\omega)}
_{T}} \|}^{2}_{L^{2}(\R)}$ is an $(L^{2})$ random variable, the Lebesgue dominated convergence theorem applies. Thus we get the  convergence $ \underset{n\to+\infty}{\lim} \int_{\Omega}  {\|{K^{(\omega)}
_{T,n}} - K^{(\omega)}_{T}\|}^{2}_{L^{2}(\R)} \ d\omega =0$; that may also be written $\ds \underset{n\to+\infty}{\lim} \int_{\R}  \E\big[ {({K
_{T,n}}(a) - K_{T}(a))}^{2}  \big] \ da =0$. This shows that there exists a strictly increasing function $\varphi:\N^{*}\rightarrow\N^{*}$ such that for almost every real $a$, $\underset{n\to+\infty}{\lim} {K_{T,\varphi(n)}}(a) = {K_{T}}(a)$, where the limit holds in $(L^{2})$. Finally, results of steps $1$ and $2$ yield, for almost every real $a$,  $K_{T}(a) = \ell^{(G)}_{T}(a)$. We have then proved that  the map $a\mapsto \ell^{(G)}_{T}(a)$  belongs to $L^{2}(\lambda \otimes \mu)$.
\vspace{0.1cm}

{\bfseries Step 3:} Proof of Equality \eqref{zpckzpkzepokzepoczepock}:

An application of the monotone class theorem allows us to establish \eqref{zpckzpkzepokzepoczepock}
only for every  $\Phi$ in $\sS(\R)$. For any $\Phi$ in  $\sS(\R)$, denote 
$\widecheck{\Phi}$ the 
function defined by setting $\widecheck{\Phi}(\xi):=\Phi(-\xi)$ 
and, for any $F$ in $\sS'(\R)$, denote $\widecheck{F}$ the 
tempered distribution defined by $<\widecheck{F},\Phi> = 
<F,\widecheck{\Phi}>$, for every $\Phi$ in $\sS(\R)$.
For every $\omega$  
in $\Omega$, define $\theta^{(\omega)}_{T}:\sS(\R)
\rightarrow\R$ by setting $\theta^{(\omega)}_{T}(\Phi):= 
\int^{T}_{0}\ \Phi (G_{s}(\omega))\ ds$. Since $
\theta^{(\omega)}_{T}$ belongs to $\sS'(\R)$, one can 
compute its Fourier transform, denoted $
\cF({\theta^{(\omega)}_{T}})$. We obtain, for almost every real $u$:
\vspace{-0.5cm}

\begin{equation}
\label{zorjozeri33}
\cF({\theta^{(\omega)}_{T}})(u) = \int^{T}_{0} e^{-iu G_{s}
(\omega)}\ ds = 2\pi\ \widecheck{f^{(\omega)}_{T}}(u)
\end{equation}

Besides, since ${\theta^{(\omega)}_{T}} = \frac{1}{2\pi} 
\cF(\cF(\widecheck{{\theta^{(\omega)}_{T}} } ))$ we get, using 
\eqref{zorjozeri33}, ${\theta^{(\omega)}_{T}} = \frac{1}{2\pi} \cF(2\pi\ 
{f^{(\omega)}_{T}}) = K^{(\omega)}_{T}$; and hence, for all $\Phi$ in $
\sS(\R)$, 
\vspace{-0.75cm}

\begin{align}
\label{zpdozeppk}
 {\theta^{(\omega)}_{T}}(\Phi) = <K^{(\omega)}_{T},\Phi>\hspace{-0.1cm}.
\end{align}

Thanks to Step 2, we know that, for almost every real $a$,  $K_{T}(a) = \ell^{(G)}_{T}(a)$. The definition of $ {\theta^{(\omega)}_{T}}$ as well as \eqref{zpdozeppk} then yield:
\vspace{-3ex}

\begin{align*}
\int^{T}_{0}\ \Phi (G_{s}(\omega))\ ds = \int_{\R} \ell^{(G)}_{T}(a)\  \Phi (a)\ da. 
\end{align*}
\vspace{-2ex}

The proof of $(ii)$ is obvious in view of the one of $(i)$, since one just has to take $f^{(\omega)}_{T}(\xi):= \frac{1}{2\pi}
\int^{T}_{0} e^{i  \xi G_{s}(\omega)} \ dR_{s}$ and $\theta^{(\omega)}_{T}(\Phi):= 
\int^{T}_{0}\ \Phi (G_{s}(\omega))\ dR_{s}$. It is therefore left to the reader.
\end{pr}

%
%


\subsection{ Comparison with the existing results on Local times of general and particular Gaussian processes}
\label{comparaison}
Let us first recall that, for any Borel set $\cB$ of $\R_{+}$, 
 the occupation measure of $G$ on $\cB$, denoted $\nu_{\cB}$, is 
 defined for all $\cD$ in $\cB(\R)$ by setting $
\nu_{\cB}(\cD):=\lambda(\{s\in \cB; \ G_{s}\in \cD \})$. When $
\nu_{\cB}$ is absolutely continuous with respect to the 
Lebesgue measure on $\R$, one defines the local time of 
$G$ on $\cB$, noted $\{L_{G}(\cB,x);\  x\in\R  \}$,  as 
being the Radon-Nikodym derivative of $\nu_{\cB}$. 
Folllowing \cite[Section $2$]{GH} we say in this case that 
$G$ is (LT). When $\cB=[0,T]$, we 
write $\{L_{G}(T,x);\  x\in\R  \}$ instead of $\{L_{G}([0,T],x);\  x\in\R  \}$.

\begin{rem}
\label{eorijeroijeoiejroej}
${\text{{\bfseries{1.}}}}$ In Theorem \ref{oedzemlijz}, instead of giving our 
proof of $(i)$, one may think to use \cite[Theorem 
($21.9$)]{GH}. In particular, we know thanks to \cite[($21.10$)]{GH} (or the equivalent condition \cite[($22.3$)]{GH}), that $G$ is (LT) and that $\{L_{G}(T,a);\  a\in\R  \}$ belongs to 
$L^{2}(\lambda \otimes \mu)$. However we do not know, 
at this stage, that $\{L_{G}(T,a);\  a\in\R  \}$ and 
$\{ \ell^{(G)}_{T}(a);\  a\in\R  \}$ are the same process (since one does not know, in particular that $\ell^{(G)}_{T}(a))$ is a random variable, which belongs to $(L^{2})$, for almost every real $a$).
Note moreover that \cite[($21.10$) \text{or} ($22.3$)]{GH} can not be used 
in the proof of $(ii)$ of Theorem \ref{oedzemlijz} either; at least without any additional 
assumption on the variance function $t\mapsto R_{t}$. In contrast, the proof we gave 
for $(i)$ can be adapted for establishing proof of $(ii)$.

${\text{{\bfseries{2.}}}}$ Note that \eqref{eovjeojoijovjeroivejro} allows one to get the occupation times formula given in \cite[Corollary p. $224$]{RY} in the case where $G$
is a continuous Gaussian semimartingale\footnote{See \cite[Proposition $5.7$]{JL17a} to get more informations on the structure of Gaussian semimartingales.}. It is also clear that \eqref{eovjeojoijovjeroivejro} generalizes, to Gaussian non-semimartingale, the occupation times formulas given in \cite[Corollary p. $224$]{RY}.

 \vspace{1ex}

 In addition to the general references cited in the 
 introduction of this section, \cite[Section 4]{LN12}, is the only one recent 
 paper in which local times of general Gaussian processes are studied 
 at such a level of generality. The authors consider therein only what we have 
 called here weighted local time, which they denote by $L_{t}(x)$ and define, when it exists, as the density of the occupation measure:
\begin{equation*}
 m_{t}(\cB) := \int^{t}_{0} \ \i1_B (G_{s})  \ dR_{s}, \hspace{2ex} \forall B\in \cB(R).
\end{equation*}
By the very definition of $L_{t}(x)$, they obtained the occupation formula \eqref{eovjeojoijovjeroivejro}. They then show that $L_{t}(x)$ is an $(L^{2})$ random variable, by giving its chaos decomposition.
    \vspace{1ex}
    
%
 In comparison with \cite[Section 4]{LN12}, and as we already stated in Section \ref{tanakaformula}, one first does not need to assume that $(t,s)\mapsto R_{t,s}$ is continuous on ${[0,T]}^{2}$, but only that $t\mapsto R_{t}$ and $t\mapsto g_{t}$ are  continuous on ${[0,T]}$.
 Besides, 
the fact that $t\mapsto R_{t}$ is of bounded variations on $[0,T]$ and such that $\alpha_{R}(\z)=0$ is 
sufficient here and one thus does not need to make Assumption (H3a) of  \cite{LN12} (\ie to assume that $t\mapsto R_{t}$ is increasing and such that $\z=\{0\}$). Finally, and as we stated in Point 3 of Remark \ref{pokerpofkerpoferpoferpofkerp}, the assumption $\E\big[\int_{\R}\hspace{-0.5ex} {\big| \hspace{-0.5ex}\int^{T}_{0}   e^{i\xi 
G_{s}}   dR_{s} \big|}^{2} d\xi 
\big]\hspace{-0.5ex}<\hspace{-0.5ex}+\infty$ is necessary and sufficient to ensure the existence of a 
weighted local time, that is an $(L^{2})$ random variable. It is then less restrictive than \cite[(H4)]{LN12}. Finally, Assumption \cite[(H4)]{LN12} is not required here but one needs to assume that $s\mapsto {R}^{\scriptscriptstyle-1/2}_{s} \in L^{1}([0,T], dR_{s})$.

%
%
%
%
%
%
%

%
\end{rem}


Let us now see
 the case of the Gaussian processes of ``reference''.

\begin{ex}
\label{oevoefijeorjiiferrefeeroivjeodzdez}
{\bfseries 1.}\ If $G$ is a Brownian motion or a Brownian bridge, then
$\ell^{(G)}_{T}(a)$, given in Definition \ref{lt}, always exists and belongs to $L^{2}(\R\times\Omega, \lambda \otimes \mu)$, as a two parameters process. Moreover, assumptions of Proposition \ref{ozrisjoerzijo}, Theorem-Definition \ref{ltw} and Theorem \ref{oedzemlijz} are all fulfilled. Of course, 
 $\ell^{(B)}_{T}$ and $\mathscr{L}^{(B)}_{T}$ both coincide and we recover classical results for Brownian motion local time (see \cite[Chapter $6$]{RY} for example) in all its generality.

{\bfseries 2.}\ {\bfseries $(i)$} If $G$ is a fBm $B^{H}$, then assumptions of Proposition \ref{ozrisjoerzijo}, of Theorem-Definition \ref{ltw} and Theorem \ref{oedzemlijz} are all fulfilled, whatever the value of 
 $H$ in $(0,1)$ is. Thus $\ell^{(G)}_{T}(a)$, given in Definition \ref{lt}, always exists and  $\{\ell^{(G)}_{T}(a)(\omega);\ (a,\omega)\}$ belongs to $L^{2}(\R\times\Omega, \lambda \otimes \mu)$.
We therefore recover the results for non weighted local time given in 
 \cite{HO02}, as well as the occupation times formulas given, for both weighted and non weighted fractional local times, in \cite[Proposition $30$ \& Definition $14$]{COU07}. In particular, if we except the occupation time formula, the results we provide here on weighted local time when $G=B^{H}$ (\ie on $\mathscr{L}^{(B^{H})}_{T}$), are new. 
 
 {\bfseries $(ii)$}  To our best knowledge, there was no proof of existence of (weighted or non-weighted) local time of general mBm (\ie in sense of  \cite[Definition 1.2]{LLVH}, recalled in our introduction) in the literature so far.
 However, several results 
exist for certain classes of mBm. More precisely, in 
\cite{BDG06,BDG07} the authors are interested in the 
mBm, which corresponds to $Y_{(1,0)}$ in \cite{StoTaq}. 
They establish, on each time interval $[a,b]\subset(0,+
\infty)$, the existence of a square integrable local time 
(\ie in $L^{2}(\R\times\Omega, \lambda \otimes \mu)$), assuming $\sup\{h(t); \ t\geq 0\} < \beta \wedge 
1$, where $h$ is $\beta$-H\"older continuous. In 
\cite{dozzi}, a Volterra-type representation of mBm is 
studied. It is established, using the sufficient condition 
given in \cite[Theorem $22.1$]
{GH}, the existence of a square integrable local time, (assuming $h$ is a differentiable and $(1/2,1)$-valued function). As 
a consequence, the authors get \eqref{zpckzpkzepokzepoczepock}. Starting from the Itô formula \cite[Theorem $2$]{nualart}, they also get 
\eqref{eovjeojoijovjeroivejro}. All these results are recovered by Theorem \ref{oedzemlijz}. In particular, the results we provide on weighted local time when $G=B^{h}$, are new.

{\bfseries 3} Finally, we provide new results on weighted and non weighted 
local times of $\mathscr{V}_{\gamma}$ - processes in the 
Proposition \ref{pofkeporfkerpofk} below.

\end{ex}

\begin{prop}
\label{pofkeporfkerpofk}
Let $B^{\gamma}$ be a $\mathscr{V}_{\gamma}$ - process. Then both 
weighted and non-weighted local times of $G$ (denoted respectively $
\ell^{(\gamma)}_{T}$ and  $\mathscr{L}^{(\gamma)}_{T}$) exist and 
belong to $L^{2}(\R\times\Omega, \lambda \otimes \mu)$, as two 
parameters processes. Moreover, for $\mu$-a.e. $\omega$ in $
\Omega$ and all $\Phi$ in $\mathscr{M}_{b}(\R)$, one has the following equalities:

\vspace{-1ex}

\begin{equation}
 \int^{T}_{0}\  \Phi(G_{s}(\omega)) \ ds =  \int_{\R}\  \ell^{(\gamma)}_{T}(y)(\omega) \ \Phi(y) \ dy.
\end{equation}
\vspace{-4ex}

\begin{equation}
\int^{T}_{0}\  \Phi(G_{s}(\omega)) \ dR_{s} =  \int_{\R}\  \mathscr{L}^{(\gamma)}_{T}(y)(\omega) \ \Phi(y) \ dy.
\end{equation}
\end{prop}

\begin{pr}
It is sufficient to apply Proposition \ref{ozrisjoerzijo}, Theorem-Definition 
\ref{ltw} and Theorem \ref{oedzemlijz}. Since $t\mapsto R_{t}$ is 
absolutely continuous (as it has been shown right above \eqref{fpkerfpekrfpo}), one 
sees, in view of Remark \ref{pokerpofkerpoferpoferpofkerp}, that one 
just has to prove that: $\E\big[\int_{\R}\ {\big| \int^{T}_{0}  \ e^{i\xi B^{\gamma}_{s}}  
\ dR_{s} \big|}^{2}\ d\xi 
\big] <+\infty$ \&  $\E\big[\int_{\R}\ {\big| \int^{T}_{0}  \ e^{i\xi B^{\gamma}_{s}}  \ 
d{s} \big|}^{2}\ d\xi \big] <+\infty$. However, thanks to \cite[(21.10) \& 
(22.3)]{GH}, one knows that this is equivalent to show that: $\int^{T}_{0} 
\int^{T}_{0} {(\Delta(t,s))}^{-1/2} \ dt \ ds<+\infty$ \&  $\int^{T}_{0} \int^{T}
_{0} {(\Delta(t,s))}^{-1/2} \ dR_{s} \ dR_{t}<+\infty$, where $\Delta(t,s):= 
\E[{(B^{\gamma}_{t}-B^{\gamma}_{s})}^{2}]$. Thanks to \cite[Proposition 
1.1]{MV05}, we know that: $R_{|t-s|} \leq \Delta(t,s) \leq 2\  R_{|t-s|}$, for 
every $(t,s)$ in $[0,T]^{2}$. We then use the same arguments as 
the one given right above \eqref{fpkerfpekrfpo}, as well as the fact that $R
$ is increasing, to get:
\vspace{-3ex}

\begin{align*}
\Theta_{T}&:=\int^{T}_{0} \int^{T}_{0} {(\Delta(s,u))}^{-1/2} \ du \ ds \leq 2\ \int^{T}_{0} \ \big(\int^{t}
_{0} {R_{u}}^{\hspace{-1ex}-1/2}   \ du\big) \ dt \leq 8 (1+R_{T}) + \frac{T^{2}}{\sqrt{\alpha}},
\end{align*}

where $\alpha:=\inf\{u\in(0,T], \ R'_{u}\leq 1 \}$ (where we set $\inf\emptyset = +\infty$). This achieves the proof.
\end{pr}

\begin{small}

\section*{Acknowledgments}
I want to express my deep gratitude to Jacques Lévy Véhel for his advices
and for the very stimulating discussions we had about this work.
I also want to thank Professor T. Hida for his warm welcome at the University of Nagoya, where a part of this paper was written, as well as Professor L. Chen and the Institute for Mathematical Sciences of Singapore (NUS), where another part of this paper was written. 

This work is dedicated to the memory of Professor Marc Yor. 


%
%
%
%
%
%
%
%
%
%
%
%

\end{small}

\vspace{3ex}


\centerline{\Large \bfseries \textsc{Appendix}}

\makeatletter
\renewcommand\theequation{\thesection.\arabic{equation}}
\@addtoreset{equation}{section}
\makeatother

\begin{appendix}

\makeatletter
\renewcommand\theequation{\thesection.\arabic{equation}}
\@addtoreset{equation}{section}
\makeatother

\vspace{-14ex}

\textcolor{white}{
\section{Appendix}
\subsection[\textcolor{black}{Bochner Integral}]{\textcolor{white}{Bochner Integral}}
\label{appendiceB}
}

\section*{\Large A.$1$  Bochner Integral}
\vspace{-1ex}

\begin{small}
The following notions about Bochner integral come from \cite[p.$72$, $80$ and $82$]{HP} and \cite[p.$247$]{Kuo2}.

\begin{defi}{Bochner integral \cite[p.$247$]{Kuo2}}
\label{bb}
Let $I$ be a Borelian subset of $\R$ endowed with the Lebesgue measure. One says that  $\Phi:I \rightarrow \ooS$ is Bochner integrable on $I$ if it satisfies the two following conditions:
\smallskip

\vspace{-1ex}

$1$  $\Phi$ is weakly measurable on $I$ \textit{i.e} $u\mapsto <\hspace{-0.2cm}<\hspace{-0.1cm} \ \Phi_{u},\varphi\  \hspace{-0.1cm}>\hspace{-0.2cm}>$  is measurable on $I$ for every $\varphi$ in $(\cS)$.
\medskip

$2$ $\exists\ p \in \bN$ such that $\Phi_{u} \in ({\cS}_{-p})$ for almost every $u \in I$ and $u \mapsto {\|\Phi_{u}\|}_{-p}$ belongs to $L^1(I)$.
\medskip

\ \ The Bochner-integral of $\Phi$ on $I$ is denoted $\int_{I}\ \Phi_{s} \ ds $ .
\end{defi}

\begin{prop}
If $\Phi:I \rightarrow \ooS$ is Bochner-integrable on $I$  then there exists an integer $p$ such that ${\left\| \int_{I} \Phi_{s}  \ ds \right\|}_{-p} \leq \int_{I} {\left\|\Phi_{s} \right\|}_{-p} \ ds$. 
Moreover $\Phi$ is also Pettis-integrable on I and both integrals coincide on I.
\end{prop}

\begin{rem}
The previous proposition shows that there is no risk of confusion by using the same notation for both
 Bochner and Pettis integrals.
\end{rem}

\begin{theo}
\label{cc}
Let $p \in \bN$ and ${(\Phi^{(n)})}_{n \in \bN}$ be a sequence of processes from $I$ to \oS such that  $\Phi^{(n)}_{u} \in ({\cS}_{-p})$ for almost every $u \in I$ and for every $n$. Assume moreover that $\Phi^{(n)}$ is  Bochner-integrable on $I$, for every $n$, and that $\lim\limits_{(n,m) \to (+\infty,+\infty)}  \int_{I}  {\big\| {\Phi}^{(m)}_{s} - {\Phi}^{(n)}_{s}\big\|}_{-p} \  ds = 0$. Then there exists an \oS-process (almost surely $({\cS}_{-p})$-valued), denoted $\Phi$, defined and Bochner-integrable on $I$, such that 

\vspace{-0.3cm}

\begin{equation}
\label{iuerhferuhfrufhu}
\lim\limits_{n \to +\infty}  \int_{I}\  {\|\Phi_{s} - {\Phi}^{(n)}_{s}\|}_{-p} \ ds = 0
\end{equation}

Furthermore, if there exists an \oS-process, denoted $\Psi$, which verifies $\eqref{iuerhferuhfrufhu}$, then $\Psi_{s} = \Phi_{s}$ for {\textit a.e.} $s$ in $I$. Finally one has  $\lim\limits_{n \to +\infty}  \int_{I}\ \Phi^{(n)}_s \ ds =  \int_{I}\ \Phi_{s} \ ds$, where the equality and the limit both hold in \oS.   
\end{theo}
\end{small}

\vspace{-8ex}

\textcolor{white}{
\subsection[\textcolor{black}{Proof of Lemmas \ref{eporkerpok} and \ref{pzokpo}}]{\textcolor{white}{Bochner Integral}}
\label{appendiceS}}

\vspace{-2ex}

\section*{\Large A.$2$  Proof of Lemmas \ref{eporkerpok} and \ref{pzokpo}}

\vspace{-1ex}


\begin{small}

{\bfseries Proof of Lemma \ref{eporkerpok}.}

%
%

The existence of the left hand side of \eqref{peoreorfeorij} results from \eqref{oeivjeorijv} and Theorem \ref{peodcpdsokcpodfckposkcdpqkoq}. Using $(ii)$ of Lemma \ref{dkdskcsdckksdksdmksdmlkskdm} and again  \eqref{oeivjeorijv}, we get, for every $\eta$ in $\sS(\R)$,
\vspace{-0.5cm}

\begin{align*}
 S(\int_{\R} \ j(y)\ \delta_{y}(<.,g>)\ dy)(\eta) &= \int_{\R} \ j(y)\ S(\delta_{y}(<.,g>)(\eta) \ dy = \frac{1}{\sqrt{2\pi} {|g|}_{0}} \int_{\R} \ j(y)\  \exp\big\{\tfrac{-{(y-<g,\eta>)}^{2}}{2 {|g|}^{2}_{0}}\big\}\ dy\\
 &= \frac{1}{\sqrt{2\pi} {|g|}_{0}} \int_{\R} \ j(y)\  \exp\big\{\tfrac{-{(y-<g,\eta>)}^{2}}{2 {|g|}^{2}_{0}}\big\}\ dy = S(j(<.,g>))(\eta),
\end{align*}

where the last equality results from \cite[Theorem $7.3$]{Kuo2}. The conclusion comes from the injectivity of $S$-transform.
\hfill $\square$
\smallskip

{\bfseries Proof of Lemma \ref{pzokpo}.}

Let $(a,\eta)$ in $\R^{*} \times \sS(\R)$ and $\varepsilon$ be in $(0,\frac{|a|}{2})$.
Among the assumptions that are required to apply the theorem of continuity 
under sign $\int$, only the domination of $\rho_{\eta}$ is not obvious and 
will therefore be stablished here. Let $
\gamma$ be chosen in $(0,\varepsilon)$ such that \eqref{epoferooerijeofoe} 
holds. Denote $\Gamma_{\gamma}:=\{u\in [0,T]; \ d(u,\Z^{T}_R)\geq 
\gamma\}$, $m_{\gamma}:=\sup\{ {R^{\scriptscriptstyle-1/2}_{u}}; \ u\in [0,T], 
\  \text{s.t. }d(u,\Z^{T}_R)\geq \gamma  \}$ and, for every $r$ in $(0,T]$, denote $M_{r}:={\sup}\{ R^{\scriptscriptstyle-1/2}_{s} \ e^{-\frac{a^{2}}{32 R_{s}}}; \ s\in [0,r] \}$. Using Lemma \ref{eofijeofjeroifjerofjer}, one gets, for every $(s,y)$ in $D(\z,
\gamma)\times D(a,\gamma)$,


%
\vspace{-3ex}

\begin{align*}
|\rho_{\eta}(y,s)| &=  |\rho_{\eta}(y,s)| \ \i1_{\Gamma_{\gamma}}(s) +  |\rho_{\eta}(y,s)| \ \i1_{D(\z,\gamma)}(s) \leq m_{\gamma} \ \i1_{\Gamma_{\gamma}}(s) + M_{\gamma}\ \i1_{D(\z,\gamma)}(s)=:f(s).\\
\end{align*}

\vspace{-4ex}

Since $f$ belongs to $L^{1}([0,T],ds)$, the domination is complete. The case $a=0$ is simpler. Indeed, one has:
\vspace{-3ex}

\begin{align*}
\forall (\varepsilon,\eta) \in \R^{*} \times \sS(\R),\hspace{0.25cm} \forall (y,s) \in D(0,\varepsilon)\times [0,T], \hspace{0.25cm}  |\rho_{\eta}(y,s)| \leq R^{\scriptscriptstyle-1/2}_{s}.
\end{align*}

Since $s\mapsto R^{\scriptscriptstyle-1/2}_{s}$  belongs to $L^{1}([0,T]\backslash\z, ds)$ by assumption, the domination is now complete and the theorem of continuity under sign $\int$ applies here and ends the proof.
\hfill $\square$
\end{small}

\vspace{-2ex}
\end{appendix}

\begin{footnotesize}

\end{footnotesize}

\end{document}